\documentclass{amsart}
\usepackage{pstricks,amssymb}

\title{Affine approach to quantum Schubert calculus}

\author{Alexander Postnikov}

\address{Department of Mathematics,
         M.I.T., Cambridge, MA 02139}

\email{apost@math.mit.edu}
\urladdr{http://www.math-mit.edu/\mytilde apost/}

\thanks{The author was supported in part by NSF grant DMS-0201494.}



\date{May~14, 2002; minor corrections on June~19, 2002}

\newtheorem{theorem}{Theorem}[section]
\newtheorem{proposition}[theorem]{Proposition}
\newtheorem{corollary}[theorem]{Corollary}
\newtheorem{definition}[theorem]{Definition}
\newtheorem{conjecture}[theorem]{Conjecture}
\newtheorem{remark}[theorem]{Remark}
\newtheorem{lemma}[theorem]{Lemma}

\def\Z{\mathbb{Z}}
\def\C{\mathbb{C}}
\def\sln{\mathfrak{sl}}
\def\slhat{\widehat{\mathfrak{sl}}}
\def\n{\mathfrak{n}}
\def\Shift{\mathrm{Shift}}
\def\Cyl{\mathcal{C}}
\def\T{\mathcal{T}}

\def\mytilde{\kern-.015in\hbox{\lower.03in\hbox{\~{}}}\kern-.01in}

\def\<{\left<}
\def\>{\right>}
\def\H{\mathrm{H}}
\def\QH{\mathrm{QH}}
\def\xx{\mathbf{x}}

\def\ll{\<\kern-2pt\<}
\def\rr{\>\kern-2pt\>}
\def\A{\mathrm{A}}
\def\e{\mathbf{e}}
\def\h{\mathbf{h}}
\def\a{\mathbf{a}}
\def\zz{\mathbf{z}}
\def\sign{\mathrm{sign}}
\def\dmin{{d_\mathrm{min}}}
\def\dmax{{d_\mathrm{max}}}
\def\Dmin{{D_\mathrm{min}}}
\def\Dmax{{D_\mathrm{max}}}
\def\da{\downarrow}
\def\diag{\mathrm{diag}}

\psset{unit=1pt}
\psset{arrowsize=4pt 1}
\psset{linewidth=1pt}
\newrgbcolor{mygreen}{.2 .7 .2}
\newrgbcolor{myred}{.7 .3 .2}
\newgray{mygray}{.90}
\countdef\x=23
\countdef\y=24
\countdef\z=25
\countdef\t=26

\def\tbox(#1,#2)#3{
\x=#1 \y=#2 
\multiply\x by 12 
\multiply\y by 12 
\z=\x \t=\y
\advance\z by 12 
\advance\t by 12 
\psline(\x,\y)(\x,\t)(\z,\t)(\z,\y)(\x,\y)
\advance\x by 6
\advance\y by 6 
\rput(\x,\y){{\bf #3}}}

\def\gbox(#1,#2)#3{
\x=#1 \y=#2 
\multiply\x by 12 
\multiply\y by 12 
\z=\x \t=\y
\advance\z by 12 
\advance\t by 12 
\psline[linecolor=gray, linewidth=0.5pt](\x,\y)(\x,\t)(\z,\t)(\z,\y)(\x,\y)
\advance\x by 6
\advance\y by 6 
\rput(\x,\y){\gray{#3}}}

\def\ebox(#1,#2)#3{
\x=#1 \y=#2 
\multiply\x by 12 
\multiply\y by 12 
\advance\x by 6
\advance\y by 6 
\rput(\x,\y){#3}}

\def\ggbox(#1,#2){
\x=#1 \y=#2 
\multiply\x by 12 
\multiply\y by 12 
\z=\x \t=\y
\advance\z by 12 
\advance\t by 12 
\psframe[fillstyle=solid, fillcolor=mygray, linewidth=0pt](\x,\y)(\z,\t)
\psline[linecolor=gray, linewidth=0.5pt](\x,\y)(\x,\t)(\z,\t)(\z,\y)(\x,\y)}

\def\tline(#1,#2)(#3,#4){
\x=#1 \y=#2 \z=#3 \t=#4
\multiply\x by 12 
\multiply\y by 12 
\multiply\z by 12 
\multiply\t by 12 
\psline(\x,\y)(\z,\t)}

\begin{document}

\begin{abstract}
This article presents a formula for products of Schubert 
classes in the quantum cohomology ring of the Grassmannian.
We introduce a generalization of Schur symmetric polynomials for shapes 
that are naturally embedded in a torus.  Then we show that the coefficients
in the expansion of these toric Schur polynomials, in terms of the regular Schur 
polynomials, are exactly the 3-point Gromov-Witten invariants; which are 
the structure constants of the quantum cohomology ring.
This construction implies that the Gromov-Witten
invariants of the Grassmannian are invariant with respect to the action of 
a twisted product of the groups $S_3$, $(\Z/n\Z)^2$, and $\Z/2\Z$.
The last group gives a certain strange duality of the quantum cohomology
that inverts the quantum parameter $q$.  
Our construction gives a solution to a problem posed by Fulton 
and Woodward about the characterization of the powers of the quantum 
parameter $q$ that occur with nonzero coefficients in the quantum product 
of two Schubert classes.  The strange duality switches the smallest 
such power of $q$ with the highest power.
We also discuss the affine nil-Temperley-Lieb algebra 
that gives a model for the quantum cohomology.
\end{abstract}

\maketitle

\tableofcontents
\section{Introduction}
\label{sec:intro}

It is well-known that the Schubert calculus is related to the theory of
symmetric functions.  The cohomology ring of the  Grassmannian is a certain
quotient of the ring of symmetric functions.  Schubert classes form a linear
basis in the cohomology and correspond to the Schur symmetric polynomials.
There is a more general class of symmetric polynomials known as the skew 
Schur polynomials.  The problem of
multiplying two Schubert classes is equivalent to the problem of expanding a
given skew Schur polynomial in the basis of ordinary Schur polynomials.  The
coefficients that appear in this expansion 
are explicitly computed using the Littlewood-Richardson rule.

Recently, in a series of papers by various authors,
attention has been drawn to the small quantum cohomology
ring of the Grassmannian.  This ring is a certain deformation of 
the usual cohomology.  Its structure constants are the 3-point
Gromov-Witten invariants, which count the numbers of 
certain rational curves of fixed degree.  

In this paper we present a quantum cohomology analogue of skew
Schur polynomials.  These are certain symmetric polynomials labelled 
by shapes that are embedded in a torus.  We show that 
the Gromov-Witten invariants are the expansion coefficients of these 
toric Schur polynomials in the basis of ordinary Schur polynomials. 

This construction implies several nontrivial results.
For example, it reproduces the known result 
that the Gromov-Witten invariants are
symmetric with respect to the action of the product of two
cyclic groups.   Also it gives a certain
``strange duality'' of the Gromov-Witten invariants that exchanges
the quantum parameter $q$ and its inverse.%
\footnote{After the original version of this paper appeared in 
the e-print arXiv, Hengelbrock informed us that he independently found 
this duality of the quantum cohomology, for $q=1$, see~\cite{HH}.}
Geometrically, this
duality implies that the number of rational curves of small degree equals 
the corresponding number of rational curves of high degree.
Another corollary of our construction is a complete characterization
of all powers of $q$ with non-zero coefficient that appear in the expansion
of the quantum product of two Schubert classes.  
This problem was posed in a recent paper by Fulton and Woodward~\cite{FW},
in which the lowest power of $q$ of was calculated.  
In virtue of strange duality, the problem of computing the highest power
of $q$ is equivalent to finding the lowest power.

The general outline of the paper follows.
In Section~\ref{sec:prelims}
we review main definitions and results related to the usual and quantum
cohomology rings of the Grassmannian.
In Section~\ref{sec:symmetric-f} we discuss symmetric functions and
their relation to the cohomology.
In Section~\ref{sec:cylindric-tableaux} we introduce our main 
tools---toric shapes and toric tableaux.  In Section~\ref{sec:Pieri-Kostka} 
we discuss the quantum Pieri formula and quantum Kostka numbers.  
In Section~\ref{sec:toric-Schur} we define toric Schur polynomials
and prove our main result on their Schur-expansion.
In Section~\ref{sec:symmetries} we discuss the cyclic symmetry and 
the strange duality of the Gromov-Witten invariants.
In Section~\ref{sec:high-low} we describe all powers of the quantum
parameter that appear in the quantum product.
In Section~\ref{sec:nTL} we discuss the action 
of the affine nil-Temperley-Lieb algebra on quantum cohomology.  
In Section~\ref{sec:final} we give final remarks,
open questions, and conjectures.

\bigskip

\textsc{Acknowledgments:}
I would like to thank Sergey Fomin, Victor Ostrik, Josh Scott, Mark
Shimozono, and Chris Woodward for interesting discussions
and helpful correspondence.  I thank Michael
Entov, whose question on quantum cohomology gave the first impulse for
writing this paper.  I am grateful to Anders Buch who pointed out on
the problem of Fulton and Woodward and made several helpful
suggestions.  I thank William Fulton for several helpful comments,
corrections, and suggestions for improvement of the exposition.

\section{Cohomology and quantum cohomology of Grassmannians}
\label{sec:prelims}

In this section we remind the reader some definitions and results related 
to cohomology and quantum cohomology rings of the Grassmannian.
An account of the classical cohomology ring of the Grassmannian
can be found in~\cite{Fulton}.
For the quantum part of the story, 
see~\cite{Agni, AW, Be, Buch, BCF, FW} and references therein.
\medskip

Let $Gr_{kn}$ be the manifold of $k$-dimensional subspaces in $\C^n$.
It is a complex projective variety called the {\it Grassmann
variety\/} or the {\it Grassmannian}.  
There is a cellular decomposition of the Grassmannian $Gr_{kn}$ into
Schubert cells $\Omega_\lambda^\circ$.  
These cells are indexed by partitions 
$\lambda$ whose Young diagrams fit inside the $k\times (n-k)$-rectangle. 
Let $P_{kn}$ be the set of such partitions.  In other words,
$$
P_{kn} = \{\lambda=(\lambda_1,\dots,\lambda_k)\mid 
n-k\geq \lambda_1\geq \dots\geq \lambda_k\geq 0\}.
$$
Recall that the Young diagram of a partition $\lambda$
is the set $\{(i,j)\in\Z^2\mid 1\leq j\leq \lambda_i\}$.  It is 
usually represented as a collection of boxes arranged on the plane 
in the same way as one would arrange elements of a matrix, see
Figure~\ref{fig:partition}. 
The boundary of the Young diagram of a partition $\lambda\in P_{kn}$ 
corresponds to a lattice path in the 
$k\times (n-k)$-rectangle from the lower left corner to the upper right corner.
Such a path can be encoded as a sequence 
$\omega(\lambda)=(\omega_1,\dots,\omega_n)$ of 0's and 1's with
$\omega_1+\cdots+\omega_n=k$, 
where 0's correspond to the right steps and 1's correspond to the upward steps
in the path, see Figure~\ref{fig:partition}.
We will say that $\omega(\lambda)$ is the {\it 01-word\/} of a partition
$\lambda\in P_{kn}$.
The 01-words are naturally associated with cosets
of the symmetric group $S_n$ modulo the maximal parabolic subgroup
$S_k\times S_{n-k}$.

\begin{figure}[ht]
\pspicture(0,0)(220,60)
\rput(-8,24){$k$}
\rput(36,58){$n-k$}
\rput(180,24){$\begin{array}{l} 
k=4,\ n=10\\
\lambda=(6,4,4,2),\ |\lambda|=16 \\
\omega(\lambda)=(0,0,1,0,0,1,1,0,0,1)\\
\lambda^\vee=(4,2,2,0),\ \lambda'=(4,4,3,3,1,1) 
 \end{array}$}
\ggbox(0,3)
\ggbox(1,3)
\ggbox(2,3)
\ggbox(3,3)
\ggbox(4,3)
\ggbox(5,3)
\ggbox(0,2)
\ggbox(1,2)
\ggbox(2,2)
\ggbox(3,2)
\ggbox(0,1)
\ggbox(1,1)
\ggbox(2,1)
\ggbox(3,1)
\ggbox(0,0)
\ggbox(1,0)
\psline[linecolor=blue, linewidth=0.5pt]{-}(0,0)(72,0)(72,48)(0,48)(0,0)
\psline[linecolor=red, linewidth=1.5pt]{-}(0,0)(24,0)(24,12)(48,12)(48,24)(48,36)(72,36)(72,48)
\psline[linecolor=red, linewidth=1.5pt]{->}(48,12)(48,24)
\pscircle*[linecolor=red](0,0){2}
\pscircle*[linecolor=red](72,48){2}
\endpspicture
\caption{A partition in $P_{kn}$}
\label{fig:partition}
\end{figure}

Let us fix a standard flag of coordinate subspaces 
$\C^1\subset\C^2\subset\cdots\subset\C^n$.
For $\lambda\in P_{kn}$ with 
$\omega(\lambda)=(\omega_1,\dots,\omega_n)$, the {\it Schubert cell\/} 
$\Omega_\lambda^\circ$ consists of all $k$-dimensional subspaces $V\subset\C^n$ 
with prescribed dimensions of intersections with the elements of the 
coordinate flag:
$\dim(V\cap \C^{i})=\omega_n+\omega_{n-1}\cdots+\omega_{n-i+1}$, 
for $i=1,\dots,k$.
The closures $\Omega_\lambda=\bar \Omega_\lambda^\circ$ 
of Schubert cells are called the {\it Schubert
varieties}.  Their fundamental cohomology classes 
$\sigma_\lambda=[\Omega_\lambda]$, $\lambda\in P_{kn}$,
called the {\it Schubert classes}, form a $\Z$-basis of the {
\it cohomology ring\/}
$\H^*(Gr_{kn})$ of the Grassmannian.  
Thus the dimension of the cohomology ring $\H^*(Gr_{kn})$ is equal to 
$|P_{kn}|=\binom nk$.
Only even-dimensional cohomologies 
of the Grassmannian may be nontrivial.  
We will use the convention that the degree of 
a cohomology class $\sigma\in \H^{2r}(Gr_{kn})$ is equal to $r$.
Then the degree of Schubert class $\sigma_\lambda$ equals 
$|\lambda|=\lambda_1+\cdots+\lambda_k$.

The basis of Schubert classes $\sigma_\lambda$ is self-dual with respect 
to the Poincar\'e pairing 
$$
\<\sigma,\rho\>=\int_{Gr_{kn}} \sigma \cdot \rho\,,\qquad
\sigma,\,\rho\in \H^*(Gr_{kn}).
$$ 
For $\lambda=(\lambda_1,\dots,\lambda_k)\in P_{kn}$, let 
$\lambda^\vee=(\lambda_1^\vee,\dots,\lambda_k^\vee)\in P_{kn}$ be the 
{\it complement partition\/} such that $\lambda_i^\vee = n-k - \lambda_{k+1-i}$,
i.e., $\lambda^\vee$ is obtained from $\lambda$ by taking the complement
to its Young diagram in the $k\times(n-k)$-rectangle and then rotating
it by $180^\mathrm{o}$ degrees, see Figure~\ref{fig:partition}.
Equivalently, if $\omega(\lambda)=(\omega_1,\dots,\omega_n)$
then $\omega(\lambda^\vee)=(\omega_n,\dots,\omega_1)$.
We have the following {\it duality theorem}
\begin{equation}
\<\sigma_\lambda,\sigma_{\mu^\vee}\>=\delta_{\lambda\mu}\qquad
\textrm{(Kronecker's delta)}.
\label{eq:lambda-mu-vee}
\end{equation}

The product of two Schubert classes $\sigma_\lambda$ and $\sigma_\mu$,
$\lambda,\, \mu\in P_{kn}$, in the cohomology ring $\H^*(Gr_{kn})$
can be written as
$$
\sigma_{\lambda}\cdot \sigma_{\mu} = \sum_{\nu}
c_{\lambda\mu}^\nu\,\, \sigma_\nu\,,
$$
where the sum is over partitions $\nu\in P_{kn}$ such that 
$|\nu|=|\lambda|+|\mu|$.
Let 
$$
c_{\lambda\mu\nu} =\int_{Gr_{kn}} \sigma_\lambda \cdot \sigma_\mu \cdot \sigma_\nu
$$
be the intersection number of the Schubert varieties
$\Omega_\lambda$, $\Omega_\mu$, and $\Omega_{\nu}$.
By~(\ref{eq:lambda-mu-vee}) we have $c_{\lambda\mu}^\nu=c_{\lambda\mu\nu^\vee}$.
This shows that the structure constants $c_{\lambda\mu}^\nu$ are nonnegative 
integer numbers.  The famous {\it Littlewood-Richardson rule\/} gives an 
explicit combinatorial 
formula for these numbers.  This is why the numbers  
$c_{\lambda\mu}^\nu=c_{\lambda\mu\nu^\vee}$ 
are usually called the {\it Littlewood-Richardson coefficients}.

\medskip

The (small) {\it quantum cohomology ring\/} $\QH^*(Gr_{kn})$ 
of the Grassmannian is 
an algebra over $\Z[q]$, where $q$ is a variable of degree $n$.  
As a linear space, the quantum cohomology is equal to the tensor product
$\H^*(Gr_{kn})\otimes \Z[q]$.  Thus Schubert classes $\sigma_\lambda$,
$\lambda\in P_{kn}$, form a $\Z[q]$-linear basis of $\QH^*(Gr_{kn})$.

The product in $\QH^*(Gr_{kn})$ is a certain $q$-deformation of the product
in $\H^*(Gr_{kn})$.  It is defined through the (3-point) Gromov-Witten
invariants.
The {\it Gromov-Witten invariant\/} $C_{\lambda\mu\nu}^d$,
usually denoted  $\<\Omega_\lambda,\Omega_\mu,\Omega_\nu\>_d$\,,
counts the number of rational curves of degree $d$ in $Gr_{kn}$ 
that meet generic translates of the Schubert varieties $\Omega_\lambda$,
$\Omega_\mu$,  and $\Omega_\nu$, provided that this number is finite.  
The last condition implies that the Gromov-Witten invariant 
$C_{\lambda\mu\nu}^d$ is defined if 
$|\lambda|+|\mu|+|\nu|=nd+k(n-k)$.  (Otherwise, we will set
$C_{\lambda\mu\nu}^d=0$.)
If $d=0$ then a degree 0 curve is a 
just a point in $Gr_{kn}$ and 
$C_{\lambda\mu\nu}^0= c_{\lambda \mu\nu}$ are the usual intersection numbers.
In general, the geometric definition of the Gromov-Witten invariants
$C_{\lambda\mu\nu}^d$ implies that they are nonnegative integer numbers.
We will use the notation $\sigma*\rho$ for the 
``quantum product'' of two classes $\sigma$ and $\rho$, i.e., their product
in the ring $\QH^*(Gr_{kn})$.
This product is a $\Z[q]$-linear operation.  
Thus it is enough to specify the quantum product
of any two Schubert classes.
It is defined as
\begin{equation}
\sigma_\lambda*\sigma_\mu = \sum_{d,\, \nu} q^d\, C_{\lambda\mu}^{\nu,d}\,
\sigma_\nu,
\label{eq:q-product}
\end{equation}
where the sum is over nonnegative integers $d$ and partitions $\nu\in P_{kn}$ 
such that $|\nu|=|\lambda|+|\mu|-d\,n$ and
the structure constants are the Gromov-Witten invariants:
$$
C_{\lambda\mu}^{\nu,d}=C_{\lambda\mu\nu^\vee}^d\,.
$$
Properties of the Gromov-Witten invariants imply that the quantum product is
a commutative and associative operation.
In the ``classical limit'' $q\to 0$, the quantum cohomology ring 
becomes the usual cohomology.  More formally,
$$
\H^*(Gr_{kn}) = \QH^*(Gr_{kn})/\<q\>.
$$

Unlike the usual Littlewood-Richardson coefficients $c_{\lambda\mu}^\nu$,
the Gromov-Witten invariants $C_{\lambda\mu}^{\nu,d}$ depend not only 
on three partitions 
$\lambda$, $\mu$, and $\nu$ but also on the numbers $k$ and $n$.  
If $n>|\lambda|+|\mu|$ then 
$C_{\lambda\mu}^{\nu,d}=\delta_{d\,0}\cdot c_{\lambda\mu}^\nu$. 
Thus all ``quantum effects'' vanish in the 
limit $n\to\infty$.

\medskip

The cohomology and quantum cohomology rings of the Grassmannian can be
presented as quotients of polynomial rings.
Let $e_1,\dots,e_k$ and $h_1,\dots,h_{n-k}$ be variables such that
degrees of $e_i$ and $h_i$ are equal to $i$.
Also let $e(t)=1+e_1 t + \cdots + e_k t^k$ and 
$h(t)=1+h_1 t + \cdots + h_{n-k} t^{n-k}$ be their generating functions.
The cohomology ring of the Grassmannian $Gr_{kn}$ can be presented 
as the following quotient of the polynomial ring:
\begin{equation}
\H^*(Gr_{kn}) = \Z[e_1,\dots,e_k,h_1,\dots,h_{n-k}]/\<e(t)\,h(-t) = 1\>,
\label{eq:H}
\end{equation}
where $\<e(t)\,h(-t) = 1\>$ is the ideal generated by the coefficients 
in the $t$-expansion of $e(t)h(-t)-1$.  More explicitly, this ideal
is generated by the expressions $e_1-h_1$, $e_2-e_1 h_1 + h_2$,
$e_3-e_2 h_1 + e_1 h_2 - h_3$, etc.  

There is a similar presentation for the quantum cohomology ring: 
\begin{equation}
\QH^*(Gr_{kn}) = \Z[q,e_1,\dots,e_k,h_1,\dots,h_{n-k}]/
\<e(t)\,h(-t) = 1+(-1)^{n-k} q\,t^n\>,
\label{eq:QH}
\end{equation}
where the ideal is generated by the coefficients 
in the $t$-expansion of the polynomial $e(t)\,h(-t)-1-(-1)^{n-k} q\,t^n$.

In this presentation of the (quantum) cohomology ring, the generator $e_i$ 
maps to the $i$-th Chern class $c_i(\mathcal{V}^*)$ of the dual to 
the universal subbundle $\mathcal{V}$ on $Gr_{kn}$ and 
$h_j$ maps to the $j$-th Chern class $c_j(\C^n/\mathcal{V})$
of the universal quotient bundle on $Gr_{kn}$.
One can always express the $h_j$ in terms of the $e_i$ 
modulo the defining ideal and vise versa.
So the (quantum) cohomology ring is generated by $e_1,\dots,e_k$
or, alternatively, by $h_1,\dots,h_{n-k}$. 
Actually, these generators are certain special Schubert classes:
$e_i = \sigma_{1^i}$ and $h_j=\sigma_j$.

The {\it Giambelli formula\/} 
shows how to express an arbitrary Schubert class $\sigma_\lambda$ 
in the cohomology ring in terms of the generators $e_i$ or $h_j$.  
According to Bertram's result~\cite{Be}, the same expression remains 
valid in the quantum cohomology ring $\QH^*(Gr_{kn})$.
Let $\lambda'=(\lambda_1',\dots,\lambda_{n-k}')\subset P_{n-k\,n}$ 
be the {\it conjugate partition\/} to $\lambda$ whose Young diagram is 
obtained by transposition of the Young diagram of $\lambda$,
see Figure~\ref{fig:partition}.
The quantum Giambelli formula claims that in the quantum cohomology
ring $\QH^*(Gr_{kn})$ we have
\begin{equation}
\sigma_\lambda = \det( h_{\lambda_{i+j-i}})_{1\leq i,\,j\leq k}
= \det( e_{\lambda_{i+j-i}'})_{1\leq i,\,j\leq n-k}\,,
\label{eq:Giambelli}
\end{equation}
where, by convention, $e_0=h_0=1$ and $e_i=h_j=0$ unless $0\leq i\leq k$
and $0\leq j\leq n-k$.
The usual Giambelli formula is obtained by setting $q=0$,
which does not change the identity due to the fact that there
are no $q$'s in it.

Let us remark that there is a natural {\it duality isomorphism\/}
of the quantum cohomology rings 
\begin{equation}
\QH^*(Gr_{kn})\simeq \QH^*(Gr_{n-k\,n}).
\label{eq:k=n-k}
\end{equation}
In this isomorphism, a Schubert class $\sigma_\lambda$ in $\QH^*(Gr_{kn})$
maps to the Schubert class $\sigma_{\lambda'}$ in $\QH^*(Gr_{n-k\,n})$
that corresponds to the conjugate partition.
In particular, the generators $e_i$ of $\QH^*(Gr_{kn})$
map to the generators $h_j$ of $\QH^*(Gr_{n-k\,n})$
and vise versa.

\section{Symmetric functions}
\label{sec:symmetric-f}

In this section we recall some facts about the ring of symmetric 
functions and its relation to the cohomology ring of the Grassmannian.
See~\cite{Mac} for more details on symmetric functions
and~\cite{Fulton} for their links with geometry. 
In the end of the section we briefly describe the approach of~\cite{BCF} 
to the quantum cohomology of the Grassmannian.
\medskip

Let $\Lambda_k=\Z[x_1,\dots,x_k]^{S_k}$ be the ring of 
{\it symmetric polynomials\/} in $x_1,\dots,x_k$.   
The ring $\Lambda$ of {\it symmetric functions\/} 
in the infinite set of variables $x_1,x_2,\dots$ is defined as
the inverse limit $\Lambda=\varprojlim \Lambda_k$
in the category of graded rings.
In other words, the elements of the ring $\Lambda$ are formal power series
(with bounded degrees) in the variables $x_1,x_2,\dots$ 
that are invariant under any finite permutation of the variables.  
The ring $\Lambda$ is freely generated by the {\it elementary 
symmetric functions\/} $e_i$ and, alternatively, by the 
{\it complete homogeneous symmetric functions\/} $h_j$:
$$
\Lambda=\Z[e_1,e_2,e_3,\dots]=\Z[h_1,h_2,h_3,\dots].
$$
These two sets of functions 
are related by the following simple equation:
$$
\left(1+\sum_{i=1}^\infty e_i\, t^i\right)\cdot
\left(1+\sum_{j=1}^\infty h_j\, (-t)^j\right) = 1.
$$
which allows one to express the $h_j$ in terms of the $e_i$ and vise versa.

For a partition $\lambda=(\lambda_1\geq\cdots\geq\lambda_l\geq0)$
and a nonnegative integer vector $\beta=(\beta_1,\dots,\beta_r)$,
a {\it semi-standard Young tableau\/} of {\it shape\/} $\lambda$ 
and {\it weight\/} $\beta$ is way to fill boxes of the 
Young diagram of shape $\lambda$ with numbers $1,\dots,r$
so that $\beta_i$ is the number of $i$'s, for $i=1,\dots,r$,
and the entries in the tableau are weakly increasing in the rows 
and strictly increasing in the columns of the Young diagram. 
For a tableau $T$ of weight $\beta$, let 
$\xx^T=\xx^\beta=x_1^{\beta_1}\cdots x_r^{\beta_r}$.

For a partition $\lambda$, 
the {\it Schur function\/} $s_\lambda$ is defined as the sum
$$
s_\lambda = s_\lambda(\xx) = \sum_{T\textrm{ of shape }\lambda} \xx^T\,
$$
over all semi-standard Young tableaux $T$
of shape $\lambda$.
The set of all Schur functions $s_\lambda$ forms a 
$\Z$-basis of the ring $\Lambda$ of symmetric functions.

The {\it Jacobi-Trudy formula\/} gives an expression of a Schur function
as the determinant of certain matrix formed by elementary symmetric
functions or, alternatively, by complete homogeneous symmetric functions.
Actually, the Jacobi-Trudy formula happened to coincide with the Giambelli 
formula for the Schubert classes.
This implies that the cohomology ring of the Grassmannian is  
isomorphic to the following quotient of the ring of symmetric functions:
\begin{equation}
\H^*(Gr_{kn}) \simeq \Lambda_{kn}=
\Lambda/\<s_\lambda\mid \lambda\not\in P_{kn}\>
=\Lambda/\<e_i, h_j\mid i>k,\, j> n-k\>\,,
\label{eq:H=Lambda}
\end{equation}
where the ideal is generated by all Schur functions whose shapes do not
fit inside the $k\times (n-k)$-rectangle.  In this isomorphism, 
the Schubert classes $\sigma_\lambda$, $\lambda\in P_{kn}$ map to the Schur 
functions $s_\lambda$.  In particular, the generators $e_i$ in~(\ref{eq:H})
map to the elementary symmetric functions $e_i(\xx)=s_{1^i}(\xx)$, and 
the generators $h_j$ in~(\ref{eq:H}) map to the complete homogeneous 
symmetric functions $h_j(\xx)=s_j(\xx)$.  
(Here and below by a slight abuse of notation we use the same letters
$e_i$ and $h_j$ for generators of $\H^*(Gr_{kn})$ and the elementary and 
complete homogeneous symmetric functions.)
This isomorphism implies that
the structure constants $c_{\lambda\mu}^\nu$ of the cohomology ring
$\H^*(Gr_{kn})$ can be 
defined in terms of the Schur functions
as follows:
$$
s_{\lambda}\cdot s_\mu = \sum_{\nu} c_{\lambda\mu}^\nu\,s_\nu\,.
$$
In other words, the Littlewood-Richardson coefficients $c_{\lambda\mu}^\nu$ are
exactly the coefficients of expansion of 
the product of two Schur functions
and, thus,
they do not depend on $k$ and $n$, provided that $\lambda,\mu,\nu\in P_{kn}$.

Suppose that $\lambda$ and $\mu$ is any pair of partitions such that
$\lambda_i\geq \mu_i$ for all $i$.  The {\it skew Young diagram\/} of shape
$\lambda/\mu$ is the set-theoretic difference of two Young diagrams of shapes
$\lambda$ and $\mu$.  As before, 
semi-standard Young tableaux of skew shape 
$\lambda/\mu$ 
and weight $\beta$
are defined as filling
of boxes of the skew Young diagram $\lambda/\mu$ with $\beta_1$ 1's,
$\beta_2$ 2's, etc.\ so that the number are weakly increasing in rows and
strictly increasing in columns; and the {\it skew Schur function\/} 
$s_{\lambda/\mu}$ is defined as the sum
$$
s_{\lambda/\mu} = s_{\lambda/\mu}(\xx) = 
\sum_{T\textrm{ of shape }\lambda/\mu} \xx^T
$$
over all semi-standard tableaux of skew shape $\lambda/\mu$.

\begin{figure}[ht]
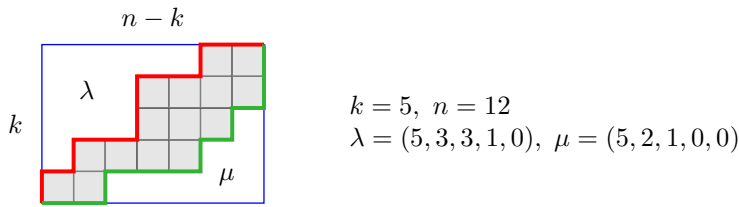

\pspicture(0,0)(230,82)
\rput(190,30){$\begin{array}{l}
k=5,\ n=12\\
\lambda=(5,3,3,1,0),\ \mu=(5,2,1,0,0) 
\end{array}
$}
\rput(-10,30){$k$}
\rput(42,70){$n-k$}
\rput(17,43){$\lambda$}
\rput(70,10){$\mu$}
\ggbox(0,0)
\ggbox(1,0)
\ggbox(1,1)
\ggbox(2,1)
\ggbox(3,1)
\ggbox(4,1)
\ggbox(3,2)
\ggbox(4,2)
\ggbox(5,2)
\ggbox(3,3)
\ggbox(4,3)
\ggbox(5,3)
\ggbox(6,3)
\ggbox(5,4)
\ggbox(6,4)
\psline[linecolor=blue, linewidth=0.5pt]{-}(0,0)(84,0)(84,60)(0,60)(0,0)
\psline[linecolor=red, linewidth=1.5pt]{-}(0,0)(0,12)(12,12)(12,24)(36,24)(36,48)(60,48)(60,60)(84,60)
\psline[linecolor=mygreen, linewidth=1.5pt]{-}(0,0)(24,0)(24,12)(60,12)(60,24)(72,24)(72,36)(84,36)(84,60)
\endpspicture
\caption{Skew shape associated with $\sigma_\lambda\cdot \sigma_\mu$}
\label{fig:skew-shape}
\end{figure}

Let $\<\cdot,\cdot\>$ the the inner product in the space of symmetric
functions $\Lambda$
such that the usual Schur functions $s_\lambda$ form an orthogonal
basis.  Then we have
$$
\<s_\lambda,\,s_\mu\cdot s_\nu\> =
\<s_{\lambda/\mu},\,s_\nu\>.
$$ 
In other words, the coefficients of expansion of a skew Schur
function in the basis of usual Schur functions are exactly the
Littlewood-Richardson coefficients:
\begin{equation}
s_{\lambda/\mu}=\sum_{\nu}c_{\mu\nu}^\lambda\,s_{\nu}\,.
\label{eq:c-lambda/mu}
\end{equation}
Equivalently, for $\lambda,\mu\in P_{kn}$, we can write
\begin{equation}
s_{\mu^\vee/\lambda}=\sum_{\nu\in P_{kn}}
c_{\lambda\mu}^{\nu}\,s_{\nu^\vee}\,.
\label{eq:mu*lambda}
\end{equation}
Here we use the fact that $c_{\lambda\mu}^{\nu}=c_{\lambda\mu\nu^\vee}$
is $S_3$-invariant under permuting of $\lambda$, $\mu$, and $\nu^\vee$.
Note that in contrast with the first formula, the second formula 
uses the complement partition operation $\nu\mapsto\nu^\vee$ and, thus, it
depends on particular values of $k$ and $n$.  The second formula means that 
the expansion coefficients of the product $\sigma_\lambda\cdot\sigma_\mu$
in the basis Schubert classes are exactly the coefficients of expansion
of the skew Schur function $s_{\mu^\vee/\lambda}$ in 
the basis of Schur functions, see Figure~\ref{fig:skew-shape}.
In this paper, we present analogues of formulas~(\ref{eq:c-lambda/mu})
and~(\ref{eq:mu*lambda}) for the quantum cohomology ring.

The Schur {\it polynomials\/} are the specializations of Schur functions
$s_\lambda(x_1,\dots,x_k)=s_\lambda(x_1,\dots,x_k,0,0,\dots)\in\Lambda_k$.
These polynomials $s_\lambda(x_1,\dots,x_k)$, where $\lambda$ ranges over 
all partitions with at most $k$ rows, form a $\Z$-basis of the ring of 
symmetric polynomials $\Lambda_k=\Z[e_1,\dots,e_k]$.
The quotient ring in~(\ref{eq:H=Lambda}) can be rewritten as
$$
\H^*(Gr_{kn}) \simeq 
\Z[e_1,\dots,e_k]/\<h_{n-k+1},\dots,h_n\>
=\Lambda_k/\<h_{n-k+1},\dots,h_n\>.
$$
The right-hand side in this equation is equivalent to
the quotient ring in~(\ref{eq:H}).

We can also present the quantum cohomology ring $\QH^*(Gr_{kn})$ 
as the following quotient of 
$\Lambda_k\otimes\Z[q]=\Z[q,e_1,\dots,e_k]$:
\begin{equation}
\QH^*(Gr_{kn})= \Lambda_k\otimes\Z[q]/I_q\,, \
\textrm{where }
I_q=\<h_{n-k+1},\dots,h_{n-1},h_n+(-1)^kq\>.
\label{eq:Gr=Lambda/ideal}
\end{equation}
This identity is equivalent to~(\ref{eq:QH}).  Notice, however, that in this 
presentation of the ring $\QH^*(Gr_{kn})$ the symmetry between the $e_i$ and 
the $h_j$ is lost.

It was shown by Bertram, Ciocan-Fontanine, and Fulton~\cite{BCF} that
the presentation $\QH^*(Gr_{kn})$ as the quotient 
ring~(\ref{eq:Gr=Lambda/ideal}) of $\Lambda_k\otimes\Z[q]$ 
allows one to write the Gromov-Witten 
invariants as alternating sums of the Littlewood-Richardson coefficients.
Indeed, in order to find the quantum product $\sigma_\lambda*\sigma_\mu$,  
$\lambda, \mu\in P_{kn}$, we need to multiply the Schur 
polynomials $s_\lambda$ and $s_\mu$ 
in the ring $\Lambda_k$; then reduce each Schur polynomial 
$s_{\tau}$ in the result modulo the ideal $I_q$; and, finally, replace
each Schur polynomial $s_{\nu}$, $\nu\in P_{kn}$, in the reduction
with the corresponding Schubert class $\sigma_{\nu}$.

An {\it $n$-rim hook\/} (a.k.a.\ border strip) is a connected skew 
Young diagram of size $n$ that contains no $2\times 2$ rectangle.  
The $n$-core of a partition $\tau$ is the partition
whose Young diagram is obtained from 
the Young diagram of $\tau$ by removing as 
many $n$-rim hooks as possible.  
It is well-known, e.g., see~\cite{Mac}, that the $n$-core 
does not depend on the order in which the rim hooks are removed.
The height of a rim hook is the number of rows it 
occupies minus 1.
For a partition $\tau$, let $d_\tau$ be the number of $n$-rim hooks 
that need to be removed from $\tau$ in order to get its $n$-core;
and let $\epsilon_\tau$ is be $d_\tau(k-1)$ minus the sum of heights 
of these rim hooks.
It is not hard to see that both $d_\tau\in\Z$ and 
$(\epsilon_\tau\ \mathrm{mod}\ 2)\in\Z/2\Z$ are well-defined.
It was shown in~\cite{BCF} 
that, for a partition $\tau$ with at most $k$ rows and $n$-core $\nu$,
the reduction of the Schur polynomial $s_{\tau}\in\Lambda_k$ modulo 
the ideal $I_q$ is equal to 
\begin{equation}
\left\{
\begin{array}{cl}
(-1)^{\epsilon_\tau}\, q^{d_\tau}\, s_{\nu} &\textrm{if }
\nu\in P_{kn},\\[.1in]
0 & \textrm{otherwise},
\end{array}
\right.
\label{eq:reduction}
\end{equation}
This claim was deduced from the Jacobi-Trudy formula. 

Thus~\cite{BCF} showed that the Gromov-Witten invariant 
$C_{\lambda\mu}^{\nu,d}$ is equal to the following alternating sum of
the Littlewood-Richardson coefficients:
$$
C_{\lambda\mu}^{\nu,d} = \sum_{\tau}(-1)^{\epsilon_\tau} 
c_{\lambda\mu}^\tau\,, 
$$
over partitions $\tau$ with at most $k$ rows 
whose $n$-core equals  $\nu$ and $d_\tau=d$.

\section{Cylindric and toric tableaux}
\label{sec:cylindric-tableaux}

In this section we develop notations needed for formulation
of our main result.  The main tool is a toric analogue of 
semi-standard Young tableaux.
\medskip

Let us fix two positive integer numbers $k$ and $n$ such that
$n>k\ge 1$ and define the {\it cylinder\/} $\Cyl_{kn}$
as the quotient 
$$
\Cyl_{kn}=\Z^2/(-k,n-k)\,\Z.
$$
In other words,
$\Cyl_{kn}$ is the quotient of the integer lattice $\Z^2$ modulo 
the action of the {\it shift operator\/} $\Shift_{kn}:\Z^2\to\Z^2$ given by 
$\Shift_{kn}:(i,j)\mapsto (i-k,j+n-k)$.  
For $(i,j)\in\Z^2$, let $\<i,j\> = (i,j)+(-k,n-k)\,\Z$ be the corresponding
element of the cylinder $\Cyl_{kn}$.
The coordinatewise partial order on $\Z^2$ gives the partial order structure
''$\preceq$'' on the cylinder $\Cyl_{kn}$ generated by the following covering relations:
$\<i,j\>\lessdot \<i,j+1\>$ and $\<i,j\>\lessdot \<i+1,j\>$.
For two points $a,b\in\Cyl_{kn}$, the {\it interval\/} $[a,b]$ is the set 
$\{c\in\Cyl_{kn}\mid a\preceq c\preceq b\}$.

\begin{definition}  A {\bf cylindric diagram} $D$ is a 
finite subset of the cylinder $\Cyl_{kn}$ closed with respect 
to the operation of taking 
intervals, i.e., for any $a,b\in D$ we have $[a,b]\subseteq D$.
\end{definition}

\begin{figure}[ht]
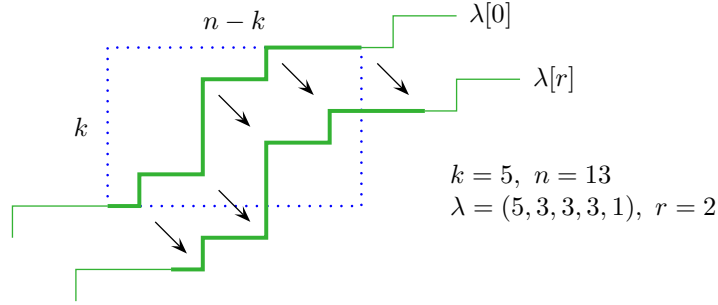

\pspicture(-12,-36)(200,77)
\rput(145,72){$\lambda[0]$}
\rput(169,48){$\lambda[r]$}
\rput(180,5){$\begin{array}{l}
k=5,\ n=13\\
\lambda=(5,3,3,3,1),\ r=2
\end{array}$}
\rput(-10,30){$k$}
\rput(48,70){$n-k$}
\psline[linecolor=blue, linewidth=1pt, linestyle=dotted]{-}(0,0)(96,0)(96,60)(0,60)(0,0)

\psline[linecolor=mygreen, linewidth=1.5pt]{-}(0,0)(12,0)(12,12)(36,12)(36,48)(60,48)(60,60)(96,60)
\psline[linecolor=mygreen, linewidth=0.5pt]{-}(96,60)(108,60)(108,72)(132,72)
\psline[linecolor=mygreen, linewidth=0.5pt]{-}(0,0)(-36,0)(-36,-12)

\psline[linecolor=mygreen, linewidth=1.5pt]{-}(24,-24)(36,-24)(36,-12)(60,-12)(60,24)(84,24)(84,36)(120,36)
\psline[linecolor=mygreen, linewidth=0.5pt]{-}(120,36)(132,36)(132,48)(156,48)
\psline[linecolor=mygreen, linewidth=0.5pt]{-}(24,-24)(-12,-24)(-12,-36)

\psline[linewidth=0.5pt]{->}(18,-6)(30,-18)
\psline[linewidth=0.5pt]{->}(42,6)(54,-6)
\psline[linewidth=0.5pt]{->}(42,42)(54,30)
\psline[linewidth=0.5pt]{->}(66,54)(78,42)
\psline[linewidth=0.5pt]{->}(102,54)(114,42)
\endpspicture
\caption{A cylindric loop $\lambda[r]$}
\label{fig:lambda[r]}
\end{figure}

Recall that a subset in a partially ordered set is 
called an {\it order ideal\/} if whenever it contains an element $a$ 
it also contains all elements which are less than $a$.
Order ideals in the partially ordered set $\Cyl_{kn}$ can be described
as follows.  
We say that an integer sequence 
$\alpha =(\dots,\alpha_{-1},\alpha_0,\alpha_1,\alpha_2,\dots)$,
infinite in both directions, 
is $(k,n)$-periodic if $\alpha_{i}=\alpha_{i+k} + (n-k)$ for any 
$i\in\Z$.
For a partition $\lambda\in P_{kn}$ and an integer number $r$,
let $\lambda[r]$ be the $(k,n)$-periodic sequence defined by
$\lambda[r]_{i+r}=\lambda_i+r$ for $i=1,\dots,k$.
All weakly decreasing $(k,n)$-periodic sequences
are of the form $\lambda[r]$. 
For any $\lambda\in P_{kn}$ and $r\in\Z$, let
$$
D_{\lambda[r]}=\{\<i,j\>\in\Cyl_{kn}\mid (i,j)\in\Z^2,\, 
j\leq \lambda[r]_i\}.
$$
The subsets $D_{\lambda[r]}$ are exactly all order ideals 
in the cylinder $\Cyl_{kn}$.
Indeed, any $\Shift_{kn}$-invariant order ideal in $\Z^2$
should be of the form $\{(i,j)\in\Z^2\mid j\leq \alpha_i\}$ 
for a weakly decreasing $(k,n)$-periodic sequence $\alpha$.
We will call the $(k,n)$-periodic sequences of the form $\lambda[r]$
{\it cylindric loops\/} of {\it type\/} $(k,n)$
because the boundary of the order ideal $D_{\lambda[r]}$ forms a closed
loop on the cylinder $\Cyl_{kn}$.
We can think about cylindric loops as infinite $\Shift_{kn}$-invariant
lattice paths on the plane.
The cylindric loop $\lambda[r]$ is obtained by shifting
the loop $\lambda[0]$ by $r$ steps to the South-East, i.e.,
by the vector $(r,r)$, see Figure~\ref{fig:lambda[r]}.

Each cylindric diagram in $\Cyl_{kn}$ is a set-theoretic difference of 
two order ideals
$$
D_{\lambda[r]/\mu[s]}=D_{\lambda[r]}\setminus D_{\mu[s]}=
\{\<i,j\>\in\Cyl_{kn}\mid (i,j)\in\Z^2,\, 
\lambda[r]_i \geq j > \mu[s]_i\},
$$
where $\lambda,\mu\in P_{kn}$ and $r,s\in\Z$. 
We will say that $D_{\lambda[r]/\mu[s]}$ is the cylindric diagram of 
{\it type\/} $(k,n)$ and {\it shape\/} $\lambda[r]/\mu[s]$.
This diagram is assumed to be empty unless 
$\lambda[r]_i\geq \mu[s]_i$ for all $i\in\Z$.  
Sometimes we will use the letter $\kappa$ 
to denote the cylindric shape $\lambda[r]/\mu[s]$
and write $D_\kappa$ instead of $D_{\lambda[r]/\mu[s]}$.
Let $|\kappa|$ denote the cardinality $|D_\kappa|$ of the cylindric diagram.

Each skew Young diagram of shape $\lambda/\mu$, with $\lambda,\mu\in P_{kn}$,
that fits inside the  $k\times (n-k)$-rectangle 
gives rise to a cylindric diagram $D_{\lambda[0]/\mu[0]}$.
In this sense we regard skew Young diagrams as a special case of cylindric
diagrams.

For two partitions $\lambda,\mu\in P_{kn}$ and a nonnegative integer $d$, 
let $\lambda/d/\mu$ be a shorthand for the cylindric shape 
$\lambda[d]/\mu[0]$. In particular, the diagrams of shape $\lambda/0/\mu$ 
are exactly the cylindric diagrams associated with a skew shape $\lambda/\mu$.
Every cylindric shape $\lambda[r]/\mu[s]$ is a shift of $\lambda/d/\mu$
by $s$ South-East steps, where $d=r-s$.  
We will often use more compact notation $\lambda/d/\mu$ 
for cylindric shapes.

Let us define rows, columns, and diagonals in the cylinder
$\Cyl_{kn}$ as follows.
The {\it $p$-th row\/} is the set $\{\<i,j\>\mid i=p\}$;
the {\it $q$-th column\/} is the set $\{\<i,j\>\mid j=q\}$;
and
the $r$-th {\it diagonal\/} 
is the set $\{\<i,j\> \mid j-i=r\}$.
The rows depend only on $p\pmod k$;
the columns depend on $q\pmod{n-k}$; and the diagonals depend 
on $r\pmod n$.
Thus the cylinder $\Cyl_{kn}$ has exactly
$k$ rows, $n-k$ columns, and $n$ diagonals.
The restriction of the partial order ``$\preceq$'' on $\Cyl_{kn}$ 
to a row, column, or diagonal gives a linear order on it.
Thus the intersection of a cylindric diagram with a row, column, or 
diagonal consists of at most one linearly ordered interval.  
These intersections 
are called rows, columns, and diagonals of the cylindric diagram.

The number of elements in the $(-k)$-th diagonal of a cylindric diagram 
$D_{\lambda[r]/\mu[s]}$ is equal to $r-s$.
In particular, the cylindric diagram of shape $\lambda/d/\mu$ has
exactly $d$ elements  in the $(-k)$-th diagonal.

\begin{definition}  
For a cylindric shape $\kappa$ and a nonnegative integer vector 
$\beta=(\beta_1,\dots,\beta_l)$, with $\beta_1+\cdots+\beta_l=|\kappa|$,
a {\bf cylindric tableau} of shape $\kappa$  and 
weight $\beta$ is a function on the cylindric diagram $T:D_\kappa\to \Z_{>0}$ 
such that $\beta_i=\#\{a\in D \mid T(a)=i\}$, for $i=1,\dots,l$. 
A cylindric tableau is {\bf semi-standard} if the function $T$ weakly
increases in the rows and strictly increases in the columns
of the cylindric diagram $D_\kappa$.
\end{definition}

Remark that {\it cylindric partitions}, which extend the notion of 
plane partitions, were introduced and studied
in~\cite{GK}.  Our semi-standard cylindric tableaux are essentially 
equivalent to {\it proper tableaux\/} from~\cite{BCF}.  
(Though the notation of~\cite{BCF} is different from ours.)

\begin{figure}[ht]
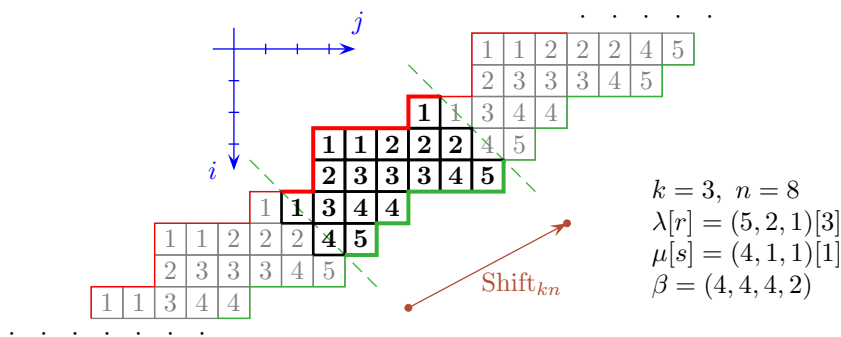

\pspicture(-48,-10)(200,120)

\rput(200,30){$\begin{array}{l}
k=3,\ n=8\\
\lambda[r]=(5,2,1)[3]\\
\mu[s]=(4,1,1)[1]\\
\beta=(4,4,4,2)
\end{array}$}

\psline[linecolor=blue, linewidth=0.5pt]{->}(-2,102)(52,102)
\psline[linecolor=blue, linewidth=0.5pt]{->}(6,110)(6,56)
\psline[linecolor=blue, linewidth=0.5pt](18,100)(18,104)
\psline[linecolor=blue, linewidth=0.5pt](30,100)(30,104)
\psline[linecolor=blue, linewidth=0.5pt](42,100)(42,104)
\psline[linecolor=blue, linewidth=0.5pt](8,90)(4,90)
\psline[linecolor=blue, linewidth=0.5pt](8,78)(4,78)
\psline[linecolor=blue, linewidth=0.5pt](8,66)(4,66)
\rput(-2,56){\blue{$i$}}
\rput(53,112){\blue{$j$}}

\psline[linecolor=mygreen, linewidth=0.5pt, linestyle=dashed]{-}(12,60)(60,12)
\psline[linecolor=mygreen, linewidth=0.5pt, linestyle=dashed]{-}(72,96)(120,48)

\psline[linecolor=myred, linewidth=0.5pt]{->}(72,4)(132,36)
\pscircle*[linecolor=myred](72,4){1.5}
\pscircle*[linecolor=myred](132,36){1.5}
\rput(115,13){\myred{$\Shift_{kn}$}}

\ebox(-1,-1){$\cdot$}
\ebox(-2,-1){$\cdot$}
\ebox(-3,-1){$\cdot$}
\ebox(-4,-1){$\cdot$}
\ebox(-5,-1){$\cdot$}
\ebox(-6,-1){$\cdot$}
\ebox(-7,-1){$\cdot$}

\gbox(-2,2){1}
\gbox(-1,2){1}
\gbox(0,2){2}
\gbox(1,2){2}
\gbox(-4,0){1}
\gbox(-3,0){1}
\gbox(-2,0){3}
\gbox(-1,0){4}
\gbox(0,0){4}
\gbox(-2,1){2}
\gbox(-1,1){3}
\gbox(0,1){3}
\gbox(1,1){3}
\gbox(2,1){4}
\gbox(1,3){1}
\gbox(2,2){2}
\gbox(3,1){5}

\gbox(8,5){4}
\gbox(9,5){5}
\gbox(7,6){1}
\gbox(8,6){3}
\gbox(9,6){4}
\gbox(10,6){4}
\gbox(8,7){2}
\gbox(9,7){3}
\gbox(10,7){3}
\gbox(11,7){3}
\gbox(12,7){4}
\gbox(13,7){5}
\gbox(8,8){1}
\gbox(9,8){1}
\gbox(10,8){2}
\gbox(11,8){2}
\gbox(12,8){2}
\gbox(13,8){4}
\gbox(14,8){5}
\ebox(11,9){$\cdot$}
\ebox(12,9){$\cdot$}
\ebox(13,9){$\cdot$}
\ebox(14,9){$\cdot$}
\ebox(15,9){$\cdot$}

\tbox(3,2){4}
\tbox(4,2){5}
\tbox(2,3){1}
\tbox(3,3){3}
\tbox(4,3){4}
\tbox(5,3){4}
\tbox(3,4){2}
\tbox(4,4){3}
\tbox(5,4){3}
\tbox(6,4){3}
\tbox(7,4){4}
\tbox(8,4){5}
\tbox(3,5){1}
\tbox(4,5){1}
\tbox(5,5){2}
\tbox(6,5){2}
\tbox(7,5){2}
\tbox(6,6){1}

\psline[linecolor=red, linewidth=0.5pt]{-}(-48,0)(-48,12)(-24,12)(-24,36)(12,36)(12,48)(24,48)
\psline[linecolor=red, linewidth=1.5pt]{-}(24,48)(36,48)(36,72)(72,72)(72,84)(84,84)
\psline[linecolor=red, linewidth=0.5pt]{-}(84,84)(96,84)(96,108)(132,108)
\psline[linecolor=mygreen, linewidth=0.5pt]{-}(0,0)(12,0)(12,12)(48,12)(48,24)
\psline[linecolor=mygreen, linewidth=1.5pt]{-}(48,24)(60,24)(60,36)(72,36)(72,48)(108,48)(108,60)
\psline[linecolor=mygreen, linewidth=0.5pt]{-}(108,60)(120,60)(120,72)(132,72)(132,84)(168,84)(168,96)(180,96)(180,108)

\endpspicture

\caption{A semi-standard cylindric tableau}
\label{fig:1}
\end{figure}

When we draw a cylindric tableau $T$ graphically, we present elements of 
its diagram as boxes and insert its values inside the boxes.  
As usual, we arrange the entries $T(\<i,j\>)$ on the plane in the same way 
as one would arrange elements of a matrix.
Figure~\ref{fig:1} gives an example of a cylindric tableau
for $k=3$ and $n=8$.
It has shape $\lambda[r]/\mu[s]=(5,2,1)[3]/(4,1,1)[1]$ 
and weight $\beta=(4,4,4,4,2)$. 
Here we presented the tableau as a $\Shift_{kn}$-periodic function 
defined on an infinite subset in $\Z^2$.
Representatives of $\Shift_{kn}$-equivalence classes of entries are
displayed in bold font.
We also indicated the $(i,j)$-coordinate system in $\Z^2$, 
the shift operator $\Shift_{kn}$, and the $(-k)$-th diagonal.
The tableau has exactly $2=3-1$ entries in this diagonal.

Let $\T_{kn}=\Z/k\Z\times \Z/(n-k)\Z$ be the integer 
{\it $k\times(n-k)$-torus}.
The torus $\T_{kn}$ is the quotient of the cylinder 
\begin{equation}
\T_{kn}=\Cyl_{kn}/(k,0)\Z=\Cyl_{kn}/(0,n-k)\Z.
\label{eq:TC}
\end{equation}
Let $\ll i,j\rr=(i+k\Z,j+(n-k)\Z)\in\T_{kn}$ be the image 
of an element $(i,j)\in\Z^2$.
Like the cylinder $\Cyl_{kn}$, the torus $\T_{kn}$ has
$k$ rows  $\{\ll i,j\rr\mid i=p\}$;
$n-k$ columns $\{\ll i,j\rr\mid j=q\}$;
and $n$ diagonals $\{\ll i,j\rr \mid j-i=r\}$.
The elements of rows, columns, and diagonals are cyclically ordered,
but there is no natural linear order on them.
Nevertheless it is still possible to define an analogue of
semi-standard tableaux whose shape is a subset of the torus $\T_{kn}$.

\begin{definition} A cylindric shape $\kappa$ is called a 
{\bf toric shape} if the restriction of the natural projection 
$p:\Cyl_{kn}\to \T_{kn}$ to the cylindric diagram of shape $\kappa$
is an injective embedding $D_\kappa\hookrightarrow\T_{kn}$.  
Two toric shapes are considered equivalent 
if they are obtained from each other by a shift by vector in $(k,0)\Z$;
thus their diagrams map to the same subset in the torus $\T_{kn}$.
A {\bf toric tableau} is a cylindric tableau of a toric shape.
\end{definition}

\begin{lemma}  A cylindric shape $\kappa$ is toric
if and only if all columns of the diagram $D_\kappa$ contain 
at most $k$ elements.  Also, a cylindric shape $\kappa$ is toric
if and only if all rows of the diagram $D_\kappa$ contain 
at most $n-k$ elements
\label{lem:toric}
\end{lemma}

\begin{proof}
Both statements immediately follow from~(\ref{eq:TC}).
\end{proof}

A cylindric loop $\lambda[r]$ 
can also be regarded as a closed loop on the torus $\T_{kn}$.
The toric shape $\lambda[r]/\mu[s]$ is formed by the elements
of the torus $\T_{kn}$ between two non-intersecting loops 
$\lambda[r]$ and $\mu[s]$.

\begin{figure}[ht]
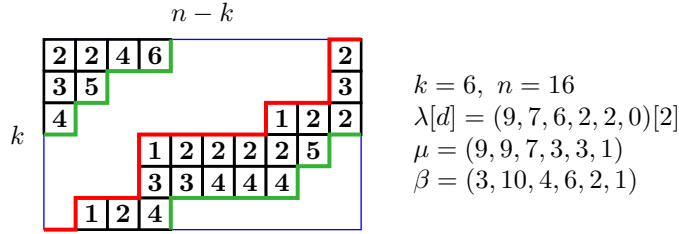

\pspicture(0,0)(230,82)

\rput(190,36){$\begin{array}{l}
k=6,\ n=16\\
\lambda[d]=(9,7,6,2,2,0)[2]\\
\mu=(9,9,7,3,3,1)\\ 
\beta=(3,10,4,6,2,1)
\end{array}$}

\rput(-10,36){$k$}
\rput(60,82){$n-k$}

\psline[linecolor=blue, linewidth=0.5pt]{-}(0,0)(120,0)(120,72)(0,72)(0,0)
\tbox(1,0){1}
\tbox(2,0){2}
\tbox(3,0){4}

\tbox(3,1){3}
\tbox(4,1){3}
\tbox(5,1){4}
\tbox(6,1){4}
\tbox(7,1){4}

\tbox(3,2){1}
\tbox(4,2){2}
\tbox(5,2){2}
\tbox(6,2){2}
\tbox(7,2){2}
\tbox(8,2){5}

\tbox(0,3){4}
\tbox(7,3){1}
\tbox(8,3){2}
\tbox(9,3){2}

\tbox(0,4){3}
\tbox(1,4){5}

\tbox(9,4){3}

\tbox(0,5){2}
\tbox(1,5){2}
\tbox(2,5){4}
\tbox(3,5){6}

\tbox(9,5){2}

\psline[linecolor=red, linewidth=1.5pt]{-}(0,0)(12,0)(12,12)(36,12)(36,36)(84,36)(84,48)(108,48)(108,72)(120,72)
\psline[linecolor=mygreen, linewidth=1.5pt]{-}(48,0)(48,12)(96,12)(96,24)(108,24)(108,36)(120,36)
\psline[linecolor=mygreen, linewidth=1.5pt]{-}(0,36)(12,36)(12,48)(24,48)(24,60)(48,60)(48,72)

\endpspicture
\caption{A semi-standard toric tableau of shape $\lambda/d/\mu$}
\label{fig:toric-tableau}
\end{figure}

The tableau given in Figure~\ref{fig:1} {\it is not\/}
a toric tableau.  It has $2$ columns with more than $3$ elements and 
two rows with more than $5$ elements.
Figure~\ref{fig:toric-tableau} gives an example of a toric tableau
drawn inside the torus $\T_{kn}$ for $k=6$ and $n=16$. 
It has shape $\lambda/d/\mu=(9,7,6,2,2,0)/2/(9,9,7,3,3,1)$ and 
weight $\beta=(3,10,4,6,2,1)$.

Notice that two different cylindric loops 
may present the same loop on the torus $\T_{kn}$.
Indeed, if we shift a loop $\lambda[r]$ by the vector $(k,0)$,
i.e., by $k$ steps to the South, we will get exactly the same
loop in the torus.

Let $\diag_0(\lambda)$ denote the number of elements in the $0$-th
diagonal of the Young diagram of a partition $\lambda\in P_{kn}$.
The number $\diag_0(\lambda)$ is also equal to the size of 
the {\it Durfee square}---the maximal square inside the Young diagram.
For a cylindric loop $\lambda[r]$, let 
$\lambda^\da[r^\da]$ be 
the cylindric loop such that $r^\da  = r+\diag_0(\lambda)$ and 
$\lambda^\da\in P_{kn}$ is the partition whose 01-word
is equal to $\omega(\lambda^\da)=(\omega_{k+1},\dots,\omega_n,
\omega_1,\dots,\omega_k)$, assuming that
$\omega(\lambda)=(\omega_1,\dots,\omega_n)$.
Notice that $\diag_0(\lambda)=\omega_{k+1}+\cdots+\omega_n$.

\begin{lemma}
For any $\lambda \in P_{kn}$ and integer $r$, two cylindric loops 
$\lambda[r]$ and $\lambda^\da[r^\da]$ present 
the same loop on the torus 
$\T_{kn}$.  Any two cylindric loops that are equivalent on the torus can be 
related by one or several transformations $\lambda[r]\mapsto 
\lambda^\da[r^\da]$.
\label{lem:loops-torus}
\end{lemma}

\begin{proof}
The cylindric loop $\lambda^\da[r^\da]$ 
is the shift of $\lambda[r]$ by the vector $(k,0)$.
\end{proof}

\section{Quantum Pieri formula and quantum Kostka numbers}
\label{sec:Pieri-Kostka}

Recall that the quantum cohomology ring $\QH^*(Gr_{kn})$ is generated
by the special Schubert classes $e_1,\dots,e_k$ or, alternatively, by 
$h_1,\dots,h_{n-k}$, see~(\ref{eq:QH}).
The quantum Pieri formula, due to Bertram~\cite{Be}, gives a rule for 
the quantum product of any Schubert class with a generator.  Thus this 
formula determines the multiplicative structure of $\QH^*(Gr_{kn})$.
In our notations we can formulate this formula as follows.

Let us say that a cylindric shape $\kappa$ is 
a {\it horizontal $r$-strip\/} if $|\kappa|=r$ and each column of
its diagram $D_\kappa$ contains at most one element.
Similarly, a cylindric shape $\kappa$ is 
a {\it vertical $r$-strip\/} if $|\kappa|=r$ and each row of $D_\kappa$ 
contains at most one element.

\begin{proposition} {\rm (Quantum Pieri formula)}
For any $\mu\in P_{kn}$ and any $i=1,\dots,k$
the quantum product $e_i*\sigma_\mu$ is given by the sum
\begin{equation}
e_i*\sigma_\mu = \sum q^d\, \sigma_\lambda
\label{eq:Pieri-e}
\end{equation}
over all $d$ and $\lambda\in P_{kn}$ such that 
$\lambda/d/\mu$ is a vertical $i$-strip.
Also, for any $j=1,\dots,n-k$,
the quantum product $h_j*\sigma_\mu$ is given by the sum
\begin{equation}
h_j*\sigma_\mu = \sum q^d\, \sigma_\lambda
\label{eq:Pieri-h}
\end{equation}
over all $d$ and $\lambda\in P_{kn}$ such that 
$\lambda/d/\mu$ is a horizontal $j$-strip.
\label{prop:Pieri}
\end{proposition}

Remark that in both cases the only possible values for $d$ are $0$ and $1$.
Indeed, any horizontal or vertical strip contains at most 1 element
in the $(-k)$-th diagonal.
Bertram proved this formula using quot schemes.
Buch gave in~\cite{Buch} a simple proof of the quantum Pieri formula using only 
the definition of Gromov-Witten invariants.
In this case we need to count degree $d=0$ curves (points)
and degree $d=1$ curves (lines) that meet generic translates of two
Schubert varieties.  This means that we need to do some linear algebra.  
For the sake of completeness we give here a short combinatorial proof 
of Proposition~\ref{prop:Pieri} in the spirit of~\cite{BCF}.

\begin{proof}
Let us first prove the rule for $e_i*\sigma_\lambda$.
Remind that $\QH^*(Gr_{kn})$ is the quotient ring~(\ref{eq:Gr=Lambda/ideal})
of $\Lambda_k\otimes\Z[q]$.   In order to find
the quantum product $e_i*\sigma_\mu$
we need to expand, in the basis of Schur polynomials,
the product $e_i\cdot s_\mu$ of the elementary symmetric polynomial with 
the Schur polynomial in the ring $\Lambda_k$;
and then reduce each Schur polynomial in the result modulo the ideal 
$I_q$ using the rule~(\ref{eq:reduction}).
The classical Pieri formula says that the product 
$e_i\cdot s_\mu$ in $\Lambda_k$ is equal to the sum
$$
e_i\cdot s_\mu = \sum s_\tau
$$
over all partitions $\tau$ with at most $k$ rows 
such that $\tau/\mu$ a vertical $i$-strip (in the classical
sense).
If $\tau\in P_{kn}$ 
then there is no $n$-rim hook that can be removed from $\tau$
and we get the terms in~(\ref{eq:Pieri-e}) with $d=0$.

Suppose that $\tau\not\in P_{kn}$. Then $\tau$ has exactly 
$n-k+1$ columns (and $\leq k$ rows). 
Thus we can remove at most one $n$-rim hook from $\tau$.
If it is impossible to remove such a rim hook then $s_\tau$ vanishes modulo
$I_q$.  Otherwise, the $n$-rim hook must occupy all $k$ rows. 
Thus, in the notation of~(\ref{eq:reduction}), 
$d_\tau=1$, $\epsilon_\tau=0$, and the $n$-core of $\tau$ is
$\lambda=(\tau_2-1,\dots,\tau_k-1,0)\in P_{kn}$.
In our notation this means that the {\it cylindric\/} shape
$\lambda/1/\mu$ is a vertical $i$-strip.
Thus we recover all terms in~(\ref{eq:Pieri-e}) with $d=1$.

The second formula~(\ref{eq:Pieri-h}) for $h_j*\sigma_\mu$
follows from the first formula~(\ref{eq:Pieri-e}) for $e_i*\sigma_\mu$
and the duality isomorphism~(\ref{eq:k=n-k}) between
$\QH^*(Gr_{kn})$ and $\QH^*(Gr_{n-k\,n})$, which switches
the $e_i$ with the $h_j$ and vertical strips with horizontal strips.
\end{proof}

Let us define the 
{\it quantum Kostka number\/} $K_{\lambda/d/\mu}^\beta$ as the number of 
semi-standard cylindric tableaux of shape $\lambda/d/\mu$ and weight $\beta$.

\begin{proposition}
For a partition $\mu\in P_{kn}$ and an integer vector 
$\beta=(\beta_1,\dots,\beta_l)$ with $0\leq \beta_i\leq n-k$,
the product $\sigma_{\mu}*h_{\beta_1}*\cdots*h_{\beta_l}$
in the quantum cohomology ring $\QH^*(Gr_{kn})$ can be expressed in terms
of the quantum Kostka numbers as follows:
$$
\sigma_{\mu}*h_{\beta_1}*\cdots*h_{\beta_l} = 
\sum_{d,\,\lambda} q^d\, K_{\lambda/d/\mu}^\beta\,\sigma_\lambda,
$$
where the sum is over nonnegative integers $d$ and partitions 
$\lambda\in P_{kn}$.
\label{prop:Kostka}
\end{proposition}

This proposition is essentially a reformulation of 
the statement on quantum Kostka numbers from~\cite[Section~3]{BCF}.
Let us show that Proposition~\ref{prop:Kostka} easily
follows from the quantum Pieri formula.

\begin{proof}
For $l=1$ the statement is equivalent to 
the quantum Pieri formula~(\ref{eq:Pieri-h}).
Indeed, cylindric tableaux of weight $\beta=(\beta_1)$ are
just horizontal $\beta_1$-strips filled with $1$'s.
Applying the quantum Pieri formula repeatedly we deduce that the coefficient
of $q^d\sigma_\lambda$ in the quantum product 
$\sigma_\mu*h_{\beta_1}*\cdots* h_{\beta_l}$
is equal to the number of chains 
$$
\lambda^{(0)}[d_0]=\mu[0],\,\lambda^{(1)}[d_1],\,\cdots,\,
\lambda^{(l-1)}[d_{l-1}],\,\lambda^{(l)}[d_{l}]=\lambda[d]
$$
such that 
$\lambda^{(i)}[d_i]/\lambda^{(i-1)}[d_{i-1}]$ is a horizontal 
$\beta_i$-strip, for $i=1,\dots,l$.
Such chains are in a one-to-one correspondence with
cylindric tableaux of shape $\lambda/d/\mu$ and weight $\beta$.
Indeed, in order to get a tableau from a chain, we just insert $i$ into 
the boxes of the $i$-th horizontal strip 
$\lambda^{(i)}[d_i]/\lambda^{(i-1)}[d_{i-1}]$, 
for $i=1,\dots,l$.
\end{proof}
 
Proposition~\ref{prop:Kostka} and the fact that
the quantum product is a commutative operation 
imply the following property of the quantum Kostka numbers.

\begin{corollary} 
The quantum Kostka numbers $K^\beta_{\lambda/d/\mu}$ 
are invariant under permuting elements $\beta_i$ of the vector $\beta$.
\label{cor:K-Sym}
\end{corollary}

It is not hard to give a direct combinatorial proof of this statement
by showing that the operators of adding horizontal (vertical) $r$-strips
to cylindric shapes commute pairwise.
The combinatorial proof is basically the same as in the case of usual 
planar shapes.

\section{Toric Schur polynomials}
\label{sec:toric-Schur}

In this section we define toric and cylindric analogues 
of Schur polynomials.  Then we formulate and prove our main result.

\medskip
For a cylindric shape $\lambda/d/\mu$,  with $\lambda,\mu\in P_{kn}$ and
$d\in\Z_{\geq0}$, we define the 
{\it cylindric Schur function\/} $s_{\lambda/d/\mu}(\xx)$ as the 
formal series in the infinite set of variables $x_1,x_2,\dots$ 
given by 
$$
s_{\lambda/d/\mu}(\xx) = 
\sum_{T} \xx^T = \sum_\beta K_{\lambda/d/\beta}\, \xx^\beta,
$$
where the first sum is over all semi-standard cylindric tableaux $T$ 
of shape $\lambda/d/\mu$; the second sum is over all possible monomials
$\xx^\beta$; and $\xx^T=\xx^\beta= x_1^{\beta_1}\cdots x_l^{\beta_l}$
for a cylindric tableau $T$ of weight $\beta=(\beta_1,\dots,\beta_l)$.

Recall that the diagrams of shape $\lambda/0/\mu$ 
are exactly the cylindric diagrams associated with a skew shape $\lambda/\mu$.
Thus
$$
s_{\lambda/0/\mu}(\xx)=s_{\lambda/\mu}(\xx)
$$ 
is the usual skew Schur function.

\begin{proposition}  
The cylindric Schur function $s_{\lambda/d/\mu}(\xx)$ belongs to the
ring $\Lambda$ of symmetric functions.
\label{prop:s-kappa-symmetric}
\end{proposition}

\begin{proof}
Follows from Corollary~\ref{cor:K-Sym}.
\end{proof}

Let us define the {\it toric Schur polynomial\/} as the specialization 
$$
s_{\lambda/d/\mu}(x_1,\dots,x_k)=s_{\lambda/d/\mu}(x_1,\dots,x_k,0,0,\dots)
$$ 
of the cylindric Schur function $s_{\lambda/d/\mu}(\xx)$.
Here, as before, $k$ is the number of rows in the torus $\T_{kn}$.
Proposition~\ref{prop:s-kappa-symmetric} implies that 
$s_{\lambda/d/\mu}(x_1,\dots,x_k)$ belongs to the ring $\Lambda_k$
of symmetric polynomials in $x_1,\dots,x_k$.
The name ``toric'' is justified by the following lemma.

\begin{lemma}
The toric Schur polynomial 
$s_{\lambda/d/\mu}(x_1,\dots,x_k)$ is nonzero if and only if
the shape $\lambda/d/\mu$ is toric.
\label{lem:schur-toric}
\end{lemma}

\begin{proof}
Let us use Lemma~\ref{lem:toric}.
If the shape $\lambda/d/\mu$ is not toric then it contains a column 
with $>k$ elements.
Thus there are no cylindric tableaux of shape $\lambda/d/\mu$ and
weight $\beta=(\beta_1,\dots,\beta_k)$ (given by a $k$-vector).  
This implies that $s_{\lambda/d/\mu}(x_1,\dots,x_k)$ is zero. 
If $\lambda/d/\mu$ is toric then all columns have $\leq k$ elements.
There are cylindric tableaux of this shape and some weight 
$\beta=(\beta_1,\dots,\beta_k)$.  For example, we can put 
the consecutive numbers $1,2,\dots$ in each column starting from the top.
This implies that $s_{\lambda/d/\mu}(x_1,\dots,x_k)\ne 0$.
\end{proof}

We are now ready to formulate our main result.
Each toric Schur polynomial $s_{\lambda/d/\mu}(x_1,\dots,x_n)$ 
can be uniquely expressed in the basis of usual Schur polynomials 
$s_\nu(x_1,\dots,x_k)$.
The next theorem links this expression to the Gromov-Witten invariants
$C_{\mu\nu}^{\lambda,d}$ that give the quantum 
product~(\ref{eq:q-product}) of Schubert classes.

\begin{theorem} 
For two partitions $\lambda,\mu\in P_{kn}$
and a nonnegative integer $d$, we have
$$
s_{\lambda/d/\mu}(x_1,\dots,x_k) = \sum_{\nu\in P_{kn}} 
C_{\mu\nu}^{\lambda,d} \,s_\nu(x_1,\dots,x_k).
$$
\label{th:main}
\end{theorem}

\begin{proof}
By the quantum Giambelli formula~(\ref{eq:Giambelli}), we have
$$
\sigma_{\mu}*\sigma_{\nu}=\sum_{w\in S_k}(-1)^{\sign(w)}\,
\sigma_{\mu}*h_{\nu_1+w_1-1}*h_{\nu_2+w_2-2}*\cdots
*h_{\nu_k+w_k-k}\,,
$$
where the sum is over all permutations $w=(w_1,\dots,w_k)$ in $S_k$.
Each of the summands in the right-hand side is given 
by Proposition~\ref{prop:Kostka}.  Extracting the coefficients of 
$q^d\sigma_\lambda$ in both sides,
we get
$$
C_{\mu\nu}^{\lambda,d}=
\sum_{w\in S_k} (-1)^{\sign(w)}\, K_{\lambda/d/\mu}^{\nu + w(\rho)-\rho},
$$
where $\nu+w(\rho)-\rho=(\nu_1+w_1-1,\dots,\nu_k+w_k-k)$.
Let us define the operator $A_\nu$ that acts on polynomials
$f\in \Z[x_1,\dots,x_k]$ as
$$
A_\nu(f) = \sum_{w\in S_k}
(-1)^{\sign(w)}\, [\textrm{coefficient of }\xx^{\nu + w(\rho)-\rho}](f).
$$
Then the previous expression can be written as
\begin{equation}
C_{\mu\nu}^{\lambda,d} = A_\nu(s_{\lambda/d/\mu}(x_1,\dots,x_k)).
\label{eq:C=A}
\end{equation}
We claim that $A_{\nu}(s_{\lambda}(x_1,\dots,x_k)) = \delta_{\lambda\nu}$.
Of course, this is a well-known identity.
This is also a special case of~(\ref{eq:C=A}) for $\mu=\emptyset$ and $d=0$.
Indeed, $C_{\emptyset\,\nu}^{\lambda,0}=c_{\emptyset\,\nu}^\lambda=
\delta_{\lambda\nu}$,
because the Schubert class $\sigma_\emptyset$ is the identity element in 
$\QH^*(Gr_{kn})$.  Thus $A_\nu(f)$ is the coefficient 
of $s_\nu$ in the expansion of $f$ in the basis of Schur polynomials.
According to~(\ref{eq:C=A}), the Gromov-Witten invariant 
$C_{\mu\nu}^{\lambda,d}$ is the coefficient
of $s_\nu$ in the expansion of $s_{\lambda/d/\mu}$, as needed.
\end{proof}

Let us reformulate our main theorem as follows.

\begin{corollary}
For two partitions $\lambda,\mu\in P_{kn}$ and a nonnegative integer $d$,
we have
$$
s_{\mu^\vee/d/\lambda}(x_1,\dots,x_k) = \sum_{\nu\in P_{kn}} 
C_{\lambda\mu}^{\nu,d}\, s_{\nu^\vee}(x_1,\dots,x_k).
$$
In other words, the coefficient of $q^d\,\sigma_{\nu}$
in the expansion of the quantum product $\sigma_\lambda*\sigma_\mu$
is exactly the same as the coefficient of $s_{\nu^\vee}$
in the Schur-expansion of the toric Schur polynomial 
$s_{\mu^\vee/d/\lambda}$.
In particular, $\sigma_\lambda*\sigma_\mu$
contains nonzero terms of the form $q^d\sigma_\nu$ if and only if
the toric Schur polynomial 
$s_{\mu^\vee/d/\lambda}$
is nonzero, i.e., $\mu^\vee/d/\lambda$
forms a valid toric shape.
\label{cor:main-variant}
\end{corollary}

\begin{proof}  The first claim is just another way to formulate
Theorem~\ref{th:main}.  Indeed, due to the $S_3$-symmetry of the Gromov-Witten
invariants, we have
$C_{\lambda\mu}^{\nu,d} =
C_{\lambda\mu\nu^\vee}^{d}=
C_{\lambda\nu^\vee}^{\mu^\vee,d}.$
The second claim follows from 
Lemma~\ref{lem:schur-toric}.
\end{proof}

This statement means that the image of the toric Schur polynomial
$s_{\mu^\vee/d/\lambda}$ in the cohomology ring $\H^*(Gr_{kn})$ under
the natural projection, see~(\ref{eq:H}), is equal to the Poincar\'e
dual of the coefficient of $q^d$ in the quantum product 
$\sigma_\lambda*\sigma_\mu$ of two Schubert classes.
In other words, the coefficient of $q^d$ in $\sigma_\lambda*\sigma_\mu$ 
is associated with the toric shape $\mu^\vee/d/\lambda$ 
in the same sense as the usual product $\sigma_\lambda\cdot \sigma_\mu$ 
is associated with the skew shape 
$\mu^\vee/\lambda$, cf.~Equation~(\ref{eq:mu*lambda}).

Theorem~\ref{th:main}
implies that all toric Schur polynomials
$s_{\lambda/d/\mu}(x_1,\dots,x_k)$ are {\it Schur-positive}, i.e., they
are positive linear combinations of usual Schur polynomials.  Indeed, the 
coefficients are the Gromov-Witten invariants, which are positive
according to their geometric definition.
Note, however, that the cylindric Schur 
function $s_{\lambda/d/\mu}(\xx)$ with non-toric shape
$\lambda/d/\mu$ may not be Schur-positive.
For example, for $k=1$ and $n=3$, we have
$$
s_{(0)/1/(0)}(\xx) = \sum_{a\leq b \leq c,\, a<c} x_a x_b x_c
= s_{21}(\xx) - s_{1^3}(\xx).
$$

\section{Symmetries of Gromov-Witten invariants}
\label{sec:symmetries}

In the section we discuss symmetries of the Gromov-Witten invariants.
We will show that they are invariant under the action
of a certain twisted product of the groups $S_3$, $(\Z/n\Z)^2$, and 
$\Z/2\Z$.  While the $S_3$-symmetry is trivial and 
the cyclic symmetry has already appeared in several papers, 
the last $\Z/2\Z$-symmetry
seems to be the most intriguing.
We call it the {\it strange duality}.
\medskip

In this section it will be convenient to use the following notation 
for the Gromov-Witten invariants: 
$$
C_{\lambda\mu\nu}(q) := 
q^d\, C_{\lambda\mu\nu}^d = q^d\, C_{\lambda\mu}^{\nu^\vee,d}\,.
$$
Recall that $d$ is determined by $\lambda$, $\mu$, and $\nu$
by $d=(|\lambda|+|\mu|+|\nu|-k(n-k))/n$.
Also let
$$
\QH^*_{\<q\>}(Gr_{kn})
=\QH^*(Gr_{kn})\otimes \Z[q,q^{-1}]
$$
be the localization of the quantum cohomology ring at
the ideal $\<q\>$.
\subsection{$S_3$-symmetry}
The invariants $C_{\lambda\mu\nu}(q)$ are symmetric 
with respect to the 6 permutations of $\lambda$, $\mu$, and $\nu$.
This is immediately clear from their geometric definition.
We have already mentioned and used this symmetry on several occasions.

\subsection{Cyclic ``hidden'' symmetry}
\label{ssec:hidden}
Let us define the {\it cyclic shift\/} operation $S$ on the set 
$P_{kn}$ of partitions as follows.  
Let $\lambda\in P_{kn}$ be a partition with the 01-word
$\omega(\lambda)=(\omega_1,\dots,\omega_n)$, see 
Section~\ref{sec:prelims}. 
Its cyclic shift $S(\lambda)$ is the partition 
$\tilde\lambda\in P_{kn}$ whose 01-word $\omega(\tilde\lambda)$ is equal to 
$(\omega_2,\omega_3,\dots,\omega_n,\omega_1)$.
Also, for the same $\lambda$,
let $\phi_i=\phi_i(\lambda)$, $i\in\Z$, be the sequence
such that $\phi_i=\omega_1+\dots+\omega_i$ for $i=1,\dots,n$;
and $\phi_{i+n}=\phi_i+k$ for any $i\in\Z$,

\begin{proposition} 
For three partitions $\lambda,\mu,\nu\in P_{kn}$ and three 
integers $a$, $b$, $c$ with $a+b+c=0$, we have
$$
C_{S^a(\lambda)\, S^b(\mu)\, S^c(\nu)}(q) =
q^{\phi_a(\lambda)+\phi_b(\mu)+\phi_c(\nu)} \, 
C_{\lambda\mu\nu}(q). 
$$
\label{prop:hidden-symmetry}
\end{proposition}
This symmetry was noticed by several people.
The first place, where it appeared in print is
Seidel's paper~\cite{Seidel}.  Agnihotri and 
Woodward~\cite[Proposition~7.2]{AW} explained
the symmetry using the Verlinde algebra.
In~\cite{P2} we gave a similar property of the quantum cohomology
of the complete flag manifold.
We call this property the {\it hidden symmetry\/} because it cannot
be detected in full generality on the level of the classical 
cohomology.  It comes from symmetries of
the extended Dynkin diagram of type $A_{n-1}$, which is an
$n$-circle.  This symmetry becomes especially transparent 
in the language of toric shapes.

\begin{proof}
It is clear from the definition that toric shapes possess
cyclic symmetry.  More precisely, for a shape $\kappa=\lambda/d/\mu$, 
the shape $S(\kappa)=S(\lambda)/\tilde d/S(\mu)$, where $\tilde d-d=
\omega_1(\mu)-\omega_1(\lambda)$, is obtained by rotation of $\kappa$.
Thus their toric Schur polynomials are the same:
$s_{\kappa}=s_{S(\kappa)}$.  
This fact, empowered by Theorem~\ref{th:main}, proves the proposition 
for $(a,b,c)=(0,1,-1)$.  The general case follows by induction from this 
claim and the $S_3$-symmetry.
\end{proof}

\begin{corollary}  For any $\lambda,\mu\in P_{kn}$ and
any integer $a$, we have the following identity
in the ring 
$\QH^*_{\<q\>}(Gr_{kn})$
$$
\sigma_{S^a(\lambda)}*\sigma_{S^{-a}(\mu)}=
q^{\phi_a(\lambda)+\phi_{-a}(\mu)}\,\sigma_\lambda*\sigma_\mu. 
$$
\label{cor:S-a-a}
\end{corollary}

The quantum cohomology ring $\QH^*(Gr_{kn})$ has the following
two ``cyclic'' Schubert classes that 
generate an $n$-dimensional algebra.
Let $E=e_k=\sigma_{1^k}$ and $H=h_{n-k}=\sigma_{n-k}$. 

\begin{proposition}
For $\lambda\in P_{kn}$ we have in the quantum cohomology ring
$$
E*\sigma_\lambda = q^{\omega_n(\lambda)} \sigma_{S^{-1}(\lambda)}
\qquad\textrm{and}\qquad
H*\sigma_\lambda = q^{1-\omega_1(\lambda)} \sigma_{S(\lambda)}.
$$
Thus the classes $E$ and $H$ generate the following 
subrings in $\QH^*(Gr_{kn})$
$$
\Z[E]/\<E^n-q^{k}\>
\qquad\textrm{and}\qquad
\Z[H]/\<H^n-q^{n-k}\>.
$$ 
The classes $E$ and $H$ are related by 
$E*H=q$.  Thus the classes $E$ and $H$ generate the same 
$\Z[q,q^{-1}]$-subalgebra in
$\QH^*_{\<q\>}(Gr_{kn})$.
Also, the class 
$$
E^{n-k}=H^k=\sigma_{(n-k)^k}
$$
is the fundamental class of a point.
\label{prop:E-H}
\end{proposition}

\begin{proof}
The first claim is a special case of the quantum Pieri formula
(Proposition~\ref{prop:Pieri}).  The remaining claims easily follow.
\end{proof}

The powers of $E$ and $H$ involve all Schubert classes $\sigma_\lambda$
with rectangular shapes $\lambda$ that have $k$ rows or $n-k$ columns.
We have 
$E^j=\sigma_{(j)^k}$ for $j=0,1,\dots,n-k$ and
$E^{n-k+i}=q^i\,\sigma_{(n-k)^{k-i}}$ for $i=0,1,\dots,k$. 
Also $H^i=\sigma_{(n-k)^i}$ for $i=0,1,\dots,k$ and
$H^{k+j}=q^j \,\sigma_{(n-k-j)^k}$ for $j=0,1,\dots,n-k$.

\subsection{Strange duality}

The quantum product has the following symmetry 
related to the Poincar\'e duality: 
$\sigma_\lambda\mapsto\sigma_{\lambda^\vee}$.

\begin{theorem}
For three partitions $\lambda,\mu,\nu\in P_{kn}$ and three 
integers $a$, $b$, $c$ with $a+b+c=n-k$, we have
$$
C_{\lambda^\vee\,\mu^\vee\,\nu^\vee}(q)= 
q^{\phi_a(\lambda)+\phi_b(\mu)+\phi_c(\nu)} \, 
C_{S^a(\lambda)\, S^b(\mu)\, S^c(\nu)}(q^{-1}).
$$
\label{thm:strange-duality}
\end{theorem}

Before we prove this theorem, let us reformulate it in
algebraic terms.
Let $D$ be the $\Z$-linear involution on 
the space $\QH^*_{\<q\>}(Gr_{kn})$ given by
$$
D:q^d\,\sigma_\lambda\longmapsto q^{-d}\,\sigma_{\lambda^\vee}.
$$
Notice that $D(1)=\sigma_{(n-k)^k}$ is the fundamental class of a point.  
It is an invertible element in the ring 
$\QH^*_{\<q\>}(Gr_{kn})$. 
By Proposition~\ref{prop:E-H}, we have $D(1)=H^k$ and
$$
(D(1))^{-1} = q^{k-n} H^{n-k}.
$$
Let us define another map $\tilde D:\QH^*_{\<q\>}(Gr_{kn})\to
\QH^*_{\<q\>}(Gr_{kn})$ as the normalization of $D$ given by
$$
\tilde D(f)= D(f)*(D(1))^{-1}.
$$
According to Proposition~\ref{prop:E-H}, the map $\tilde D$ 
is explicitly given by
$$
\tilde D:q^d\,\sigma_\lambda \longmapsto  q^{-d - \diag_0(\lambda)}\, 
\sigma_{S^{n-k}(\lambda^\vee)}\,,
$$
where $\diag_0(\lambda)=\phi_{n-k}(\lambda^\vee)=k-\phi_k(\lambda)$ 
is the size of the $0$-th diagonal of the Young diagram of $\lambda$.

\begin{theorem}
The map $\tilde D$ is a homomorphism of the ring 
$\QH^*_{\<q\>}(Gr_{kn})$.  The map $\tilde D$ is also an involution.
It inverts the quantum parameter: $\tilde D(q) = q^{-1}$.
\label{th:homomorph}
\end{theorem}

For $q=1$, the involution $\tilde D$ was independently discovered 
from a different point of view by Hengelbrock~\cite{HH}.  He showed 
that it comes from complex conjugation of the points in 
$\mathrm{Spec}\,R$, where $R=\QH^*(Gr_{kn})/\<q-1\>$.

The claim that $\tilde D$ is an involution of $\QH^*_{\<q\>}(Gr_{kn})$
implies that the map $\lambda\mapsto \tilde\lambda=S^{n-k}(\lambda^\vee))$
is an involution on partitions in $P_{kn}$ and  
$\diag_0(\lambda)=\diag_0(\tilde \lambda)$.
It is easy to see this combinatorially.
Indeed, if the 01-word of $\lambda$ is 
$\omega(\lambda)=(\omega_1,\dots,\omega_n)$
then the 01-word of $\tilde\lambda$ is 
$\omega(\tilde \lambda)=
(\omega_k,\omega_{k-1},\dots,
\omega_1,\omega_n,\omega_{n-1},\dots,\omega_{k+1})$
and $\diag_0(\lambda)=\diag_0(\tilde\lambda)=\omega_{k+1}+\dots+\omega_n$.

The previous theorem is equivalent to the following property of
the involution~$D$.

\begin{proposition} 
We have the identity
$$
D(f*g)*D(h)=D(f)*D(g*h),
$$
for any $f, g, h\in\QH^*_{\<q\>}(Gr_{kn})$. 
\label{prop:D(f*g)}
\end{proposition}

We will need the following lemma.

\begin{lemma}
For any $f\in\QH^*_{\<q\>}(Gr_{kn})$ and any 
$i=0,\dots,k$, we have
\begin{equation}
D(f*e_i)=q^{-1}\,D(f)*h_{n-k}*e_{k-i}\,.
\label{ea:Dfei}
\end{equation}
Here we assume that $e_0=1$.
\label{lem:Dfei}
\end{lemma}

\begin{proof}  
Since $D(q^k\, f)=q^{-k}\, D(f)$, it is enough to prove the lemma
for a Schubert class $f=\sigma_\lambda$. 
According to the quantum Pieri formula (Proposition~\ref{prop:Pieri}),
$\sigma_\lambda*e_i$ is given by the sum over all possible ways 
to add a vertical $i$-strip to the cylindric loop $\lambda[0]$.
Thus $D(\sigma_\lambda*e_i)$ is given by the sum over all possible ways 
to {\it remove\/} a vertical $i$-strip from $\lambda^\vee[0]$.  
In other words, we have
$$
D(\sigma_\lambda*e_i)=\sum q^{-d}\, \lambda_\mu\,,
$$
where the sum is 
over $\mu\in P_{kn}$ and $d$ such that $\lambda^\vee/d/\mu$ is a vertical
$i$-strip.
By Proposition~\ref{prop:E-H}, the right-hand side 
of~(\ref{ea:Dfei}) is equal to
$$
q^{-1}\,\sigma_{\lambda^\vee}*h_{n-k}*e_{k-i}=
q^{-\omega_1(\lambda^\vee)}\,\sigma_{S(\lambda^\vee)}*e_{k-i}\,.
$$
We obtain exactly the same expressions in both cases.
Indeed, {\it removing\/} a vertical $i$-strip from a cylindric shape 
means exactly the same as {\it cyclically shifting\/} the shape and 
then {\it adding\/} a vertical $(k-i)$-strip.  
By looking on the formula for a minute, we also see that the powers 
of $q$ in both cases are equal to each other.
\end{proof}

\begin{proof}[Proof of Proposition~\ref{prop:D(f*g)}]
Again, since multiplying $g$ by a power of $q$ does not change the formula,
it is enough to verify the statement when $g$ belongs to some set that
spans the algebra $\QH^*_{\<q\>}(Gr_{kn})$ over $\Z[q,q^{-1}]$.  Let us prove 
the statement when $g=e_{i_1}*e_{i_2}*\cdots*e_{i_l}$.  If $l=1$,
then by Lemma~\ref{lem:Dfei}, we have
$$
D(f*e_i)*D(h)=q^{-1}\,h_{n-k}*e_{k-1}*D(f)*D(h)=D(f)*D(e_i*h).
$$
The general case follows from this case.  We just need to move the $l$
factors $e_{i_1}, \dots, e_{i_l}$ one by one from 
the first $D$ to the second $D$.
\end{proof}

\begin{proof}[Proof of Theorem~\ref{th:homomorph}]
Proposition~\ref{prop:D(f*g)} with $h=1$ says that
$D(f*g)*D(1)=D(f)*D(g)$.  It is equivalent to saying that the
normalization $\tilde D$ is a homomorphism. 
We already proved combinatorially that $\tilde D$ is an involution.
Let us also it deduce this fact algebraically
from Proposition~\ref{prop:D(f*g)}:
$$
\tilde D(\tilde D(f))= \frac{D\left(\frac{D(f)}{D(1)}\right)*D(D(1))}{D(1)}=
\frac{D(D(f))*D\left(\frac{D(1)}{D(1)}\right)}{D(1)} = f.
$$
The fact that $\tilde D(q)=q^{-1}$ is clear from the definition.
\end{proof}

\begin{corollary}
The coefficient of $q^d\,\sigma_{\nu^\vee}$ in the quantum
product $\sigma_{\lambda^\vee}*\sigma_{\mu^\vee}$
is equal to the coefficient of $q^{\diag_0(\nu)-d}\,\sigma_{S^k(\nu)}$ 
in the quantum product $\sigma_{\lambda}*\sigma_{\mu}$.
\label{cor:duality-vee*vee}
\end{corollary}

\begin{proof}
By setting $f=\sigma_{\lambda^\vee}$, $g=\sigma_{\mu^\vee}$, and $h=1$
in Proposition~\ref{prop:D(f*g)}, we obtain
$$
D(\sigma_{\lambda^\vee}*\sigma_{\mu^\vee})*D(1)=\sigma_\lambda*\sigma_\mu.
$$
Since $D(1)=\sigma_{(n-k)^k}$ is the fundamental class of a point,
we get, by Proposition~\ref{prop:E-H}, 
$$
D(q^d*\sigma_{\nu^\vee})*D(1)=q^{-d}*\sigma_\nu*\sigma_{(n-k)^k}=
q^{-d}\, H^k*\sigma_{\nu}= q^{\diag_0(\nu)-d} \,\sigma_{S^k(\nu)}\,.
$$
Here we used the fact that $\diag_0(\nu)=k-\phi_k(\nu)$.
\end{proof}

We can now prove the first claim of this subsection.

\begin{proof}[Proof of Theorem~\ref{thm:strange-duality}]
Corollary~\ref{cor:duality-vee*vee} is equivalent 
to the special case of Theorem~\ref{thm:strange-duality} for
$a=b=0$ and $c=n-k$. 
The general case follows by  Proposition~\ref{prop:hidden-symmetry}.
\end{proof}

The statement of Corollary~\ref{cor:duality-vee*vee} means 
that the terms in the expansion
of the quantum product $\sigma_\lambda*\sigma_\mu$ 
are in one-to-one correspondence
with the terms in the expansion of the quantum product 
$\sigma_{\lambda^\vee}*\sigma_{\mu^\vee}$
so that the coefficients of corresponding terms are equal to each other.
Notice that terms with low powers of $q$ correspond to terms with 
high powers of $q$ and vise versa.
This property seems mysterious from the point of view of quantum cohomology.
Why should the number of some rational curves of high degree be equal 
to the number of some rational curves of low degree?
This strange duality is also hidden on the classical level.
Indeed, if $|\lambda|+|\mu|<\dim_\C Gr_{kn}$, then 
the product $\sigma_\lambda\cdot \sigma_\mu$ is always nonzero 
and  the product $\sigma_{\lambda^\vee}\cdot\sigma_{\mu^\vee}$ 
always vanishes in the classical cohomology ring $\H^*(Gr_{kn})$.

\medskip
Let us reformulate this duality in terms of toric Schur polynomials.
For a toric shape $\kappa=\lambda[r]/\mu[s]$, let us define
the {\it complement toric shape\/} as
$$
\kappa^\vee=\mu^\da[s^\da]/\lambda[r], 
$$
where the transformation $\mu[s]\mapsto \mu^\da[s^\da]$ 
is the same as in Lemma~\ref{lem:loops-torus}.

\begin{figure}[ht]
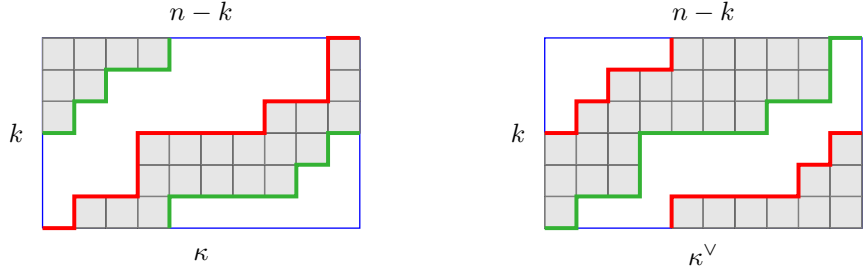

\pspicture(-10,-10)(120,82)
\rput(-10,36){$k$}
\rput(60,82){$n-k$}
\rput(60,-10){$\kappa$}
\psline[linecolor=blue, linewidth=0.5pt]{-}(0,0)(120,0)(120,72)(0,72)(0,0)
\ggbox(1,0)
\ggbox(2,0)
\ggbox(3,0)
\ggbox(3,1)
\ggbox(4,1)
\ggbox(5,1)
\ggbox(6,1)
\ggbox(7,1)
\ggbox(3,2)
\ggbox(4,2)
\ggbox(5,2)
\ggbox(6,2)
\ggbox(7,2)
\ggbox(8,2)
\ggbox(0,3)
\ggbox(7,3)
\ggbox(8,3)
\ggbox(9,3)
\ggbox(0,4)
\ggbox(1,4)
\ggbox(9,4)
\ggbox(0,5)
\ggbox(1,5)
\ggbox(2,5)
\ggbox(3,5)
\ggbox(9,5)
\psline[linecolor=red, linewidth=1.5pt]{-}(0,0)(12,0)(12,12)(36,12)(36,36)(84,36)(84,48)(108,48)(108,72)(120,72)
\psline[linecolor=mygreen, linewidth=1.5pt]{-}(48,0)(48,12)(96,12)(96,24)(108,24)(108,36)(120,36)
\psline[linecolor=mygreen, linewidth=1.5pt]{-}(0,36)(12,36)(12,48)(24,48)(24,60)(48,60)(48,72)
\endpspicture
\qquad \qquad \qquad
\pspicture(-10,-10)(120,82)
\rput(-10,36){$k$}
\rput(60,82){$n-k$}
\rput(60,-10){$\kappa^\vee$}
\psline[linecolor=blue, linewidth=0.5pt]{-}(0,0)(120,0)(120,72)(0,72)(0,0)
\ggbox(0,0)
\ggbox(0,1)
\ggbox(1,1)
\ggbox(2,1)
\ggbox(0,2)
\ggbox(1,2)
\ggbox(2,2)
\ggbox(1,3)
\ggbox(2,3)
\ggbox(3,3)
\ggbox(4,3)
\ggbox(5,3)
\ggbox(6,3)
\ggbox(2,4)
\ggbox(3,4)
\ggbox(4,4)
\ggbox(5,4)
\ggbox(6,4)
\ggbox(7,4)
\ggbox(8,4)
\ggbox(4,5)
\ggbox(5,5)
\ggbox(6,5)
\ggbox(7,5)
\ggbox(8,5)
\ggbox(4,0)
\ggbox(5,0)
\ggbox(6,0)
\ggbox(7,0)
\ggbox(8,0)
\ggbox(9,0)
\ggbox(8,1)
\ggbox(9,1)
\ggbox(9,2)
\psline[linecolor=mygreen, linewidth=1.5pt]{-}(0,0)(12,0)(12,12)(36,12)(36,36)(84,36)(84,48)(108,48)(108,72)(120,72)
\psline[linecolor=red, linewidth=1.5pt]{-}(48,0)(48,12)(96,12)(96,24)(108,24)(108,36)(120,36)
\psline[linecolor=red, linewidth=1.5pt]{-}(0,36)(12,36)(12,48)(24,48)(24,60)(48,60)(48,72)
\endpspicture
\caption{The complement toric shape}
\label{fig:toric-complement}
\end{figure}

This definition has the following simple geometric meaning.
The image of the diagram $D_{\kappa^\vee}$ of shape $\kappa^\vee$ 
in the torus $\T_{kn}$ is the complement 
to the image of the diagram $D_\kappa$ of shape $\kappa$, see 
Figure~\ref{fig:toric-complement}.
If $\kappa=\lambda/d/\mu$ then $\kappa^\vee$ obtained by a shift of 
the toric shape $\mu^\da/d'/\lambda$, where
$\mu^\da=S^k(\mu)$ and $d'=\diag_0(\mu)-d=
\phi_{n-k}(\mu^\vee)-d$.  Thus the toric Schur polynomial
$s_{(\lambda/d/\mu)^\vee}$ is equal to 
$s_{\mu^\da/d'/\lambda}$.

\begin{corollary}
For any toric shape $\kappa$, 
the coefficients in the Schur-expansion of 
the toric Schur polynomial $s_\kappa$ correspond to the coefficients
in the Schur-expansion of the toric Schur polynomial $s_{\kappa^\vee}$,
as follows:
$$
s_\kappa=\sum_{\nu\in P_{kn}} a_\nu\, s_\nu
\qquad\textrm{has the same coefficients $a_\nu$ as in}\qquad
s_{\kappa^\vee}=\sum_{\nu\in P_{kn}} a_\nu\, s_{\nu^\vee}\,.\\[.1in]
$$
\end{corollary}

\begin{proof}
Suppose that $\kappa=\lambda/d/\mu$.
By Theorem~\ref{th:main}, 
the coefficient of $s_\nu$ in the Schur-expansion of 
$s_\kappa$ is equal to $C_{\mu\nu}^{\lambda,d}=C_{\mu\nu\lambda^\vee}^d$.
On the other hand, the coefficient of $s_{\nu^\vee}$ in
the Schur-expansion of $s_{\kappa^\vee}=s_{\mu^\da/d'/\lambda}$
is equal to 
$C_{\lambda\,\nu^\vee}^{\mu^\da,d'}=
C_{\lambda\,\nu^\vee S^{n-k}(\mu^\vee)}^{d'}$.
The equality of these two coefficients is a special case of 
Theorem~\ref{thm:strange-duality}.
\end{proof}

\subsection{Essential interval}
\label{ssec:essential}

In many cases the hidden symmetry and the strange duality
imply that Gromov-Witten invariant vanishes.  In some other cases
these symmetries allow us to reduce Gromov-Witten invariant to 
Littlewood-Richardson coefficient. 
For three partitions $\lambda,\mu,\nu\in P_{kn}$, let us 
define three numbers 
$$
\begin{array}{l}
\displaystyle
\dmin(\lambda,\mu,\nu)=-
\min_{a+b+c=0} (\phi_a(\lambda)+\phi_{b}(\mu)+\phi_c(\nu)),
\\[.15in]
\displaystyle
\dmax(\lambda,\mu,\nu)=
-\max_{a+b+c=k-n} (\phi_a(\lambda)+\phi_{b}(\mu)+\phi_c(\nu)),\\[.15in]
d(\lambda,\mu,\nu)=(|\lambda|+|\mu|+|\nu|-k(n-k))/n,
\end{array}
$$
where in the first and in the second cases 
the maximum and minimum is taken over all
triples of integers $a$, $b$, and $c$ that satisfy the given condition.

\begin{proposition}
Let $\lambda,\mu,\nu\in P_{kn}$ be three partitions and
$\dmin=\dmin(\lambda,\mu,\nu)$, 
$\dmax=\dmax(\lambda,\mu,\nu)$.
Then the Gromov-Witten invariant $C_{\lambda\mu\nu}^d$ is equal to zero
unless $d=d(\lambda,\mu,\nu)$ and $\dmin\leq d \leq \dmax$.  
If $d=\dmin$ and $(a,b,c)$ 
is a triple such that $a+b+c=0$ and
$d=-(\phi_a(\lambda)+\phi_b(\mu)+
\phi_c(\nu))$ then
$$
C_{\lambda\mu\nu}^\dmin= c_{S^a(\lambda) S^b(\mu) S^c(\nu)}.
$$
Similarly, if $d=\dmax$ and $(a,b,c)$ is a triple such that 
$a+b+c=k-n$ and 
$d=-(\phi_a(\lambda)+\phi_b(\mu)+
\phi_c(\nu))$ then
$$
C_{\lambda\mu\nu}^\dmax= c_{S^{-a}(\lambda^\vee) S^{-b}(\mu^\vee) 
S^{-c}(\nu^\vee)}.
$$
\label{prop:essential}
\end{proposition}

\begin{proof}
The claim that $C_{\lambda\mu\nu}^d=0$ unless $d=d(\lambda,\mu,\nu)$ 
follows directly from the definition of the Gromov-Witten invariants.
Proposition~\ref{prop:hidden-symmetry} says that
$C_{\lambda\mu\nu}^d = 
C_{S^a(\lambda) S^b(\mu) S^c(\nu)}^{\tilde d}$,
where $\tilde d=d+\phi_a(\lambda)+\phi_b(\mu)+\phi_c(\nu)$.
The Gromov-Witten invariant in the right-hand side vanishes if $\tilde d<0$ 
and it is Littlewood-Richardson coefficient if $\tilde d=0$.  
This proves that $C_{\lambda\mu\nu}^d=0$ for $d< \dmin$ and 
that $C_{\lambda\mu\nu}^\dmin$ is Littlewood-Richardson coefficient. 
Similarly, the statement that $C_{\lambda\mu\nu}^d=0$ for 
$d>\dmax$ and $C_{\lambda\mu\nu}^\dmax$ is 
Littlewood-Richardson coefficient is a consequence of 
Theorem~\ref{thm:strange-duality}.
\end{proof}

We say that the integer interval 
$[\dmin(\lambda,\mu,\nu),\dmax(\lambda,\mu,\nu)]$ is
the {\it essential interval\/} for the triple of partitions
$\lambda,\mu,\nu\in P_{kn}$.

\section{Powers of $q$ in quantum product of Schubert classes}
\label{sec:high-low}

In this section we discuss the following problem:  
what is the set of all powers $q^d$ that appear with non-zero coefficients
in the Schubert-expansion of a given quantum product 
$\sigma_\lambda*\sigma_\mu$ ?
The lowest such power of $q$ was established in~\cite{FW}.  
Some bounds for the highest power of $q$ were found in~\cite{Yong}.
In this section we present a simple answer to this problem.  We would like 
to thank here Anders Buch who remarked that our main theorem resolves this 
problem and made several helpful suggestions.

\medskip
We have already formulated the answer to this problem in 
Corollary~\ref{cor:main-variant}.  The quantum product 
$\sigma_\lambda*\sigma_\mu$ contains nonzero terms with given power
$q^d$ if and only if $\mu^\vee/d/\lambda$ forms a valid toric shape.
Let us spell out this last condition explicitly.

Recall that, 
for $\lambda\in P_{kn}$ with 01-word 
$\omega(\lambda)=(\omega_1,\dots,\omega_n)$,
the sequence $\phi_i(\lambda)$, $i\in\Z$, is defined by 
$\phi_i(\lambda)=\omega_1+\dots+\omega_i$ for $i=1,\dots,n$
and $\phi_{n+i}(\lambda)=\phi_i(\lambda)+k$ for any $i\in\Z$,
see Section~\ref{sec:symmetries}.
For any $\lambda,\mu\in P_{kn}$, 
let us define two numbers $\Dmin$ and $\Dmax$ as follows:
$$
\begin{array}{l}
\displaystyle
\Dmin=
\Dmin(\lambda,\mu)=-
\min_{i+j=0} (\phi_i(\lambda)+\phi_{j}(\mu)),
\\[.15in]
\displaystyle
\Dmax=
\Dmax(\lambda,\mu)=
-\max_{i+j=k-n} (\phi_i(\lambda)+\phi_{j}(\mu)),
\end{array}
$$
where in the both cases the maximum or minimum is taken over all
integers $i$ and $j$ that satisfy the given condition, 
cf.~Section~\ref{ssec:essential}.

\begin{theorem}
For any pair $\lambda,\mu\in P_{kn}$, we have
$\Dmin\leq \Dmax$ and the set of all $d$'s
such that the power $q^d$ appears in $\sigma_\lambda*\sigma_\mu$
with nonzero coefficient is exactly the integer interval
$\Dmin\leq d \leq \Dmax$.
In particular, the quantum product  $\sigma_\lambda*\sigma_\mu$
is always nonzero.
\label{th:dmin-max}
\end{theorem}

The claim that $\Dmin$ is the lowest power of
$q$ with non-zero coefficient is due to Fulton and Woodward~\cite{FW}.
Some bounds for the highest power of $q$ 
were given by Yong in~\cite{Yong}.  He also formulated a conjecture
that the powers of $q$ that appear in the expansion of the quantum 
product $\sigma_\lambda*\sigma_\mu$ form an interval of consecutive
integers.

The number $\Dmin$ was defined in~\cite{FW} in terms
of overlapping diagonals in two Young diagrams.  
Let us show how to reformulate
our definitions of $\Dmin$ and $\Dmax$ in these terms.
For a partition $\lambda\in P_{kn}$ and $i=-k,\dots,n-k$,  let 
$\diag_i(\lambda)$
denote the number of elements in the $i$-th diagonal of the Young 
diagram of shape $\lambda$.  In particular, 
$\diag_{-k}(\lambda)= \diag_{n-k}(\lambda)=0$.
Then, for any pair of partitions $\lambda,\mu\in P_{kn}$, 
we have
$$
\begin{array}{l}
\displaystyle
\Dmin(\lambda,\mu)=\max_{i=-k,\dots,n-k} (\diag_i(\lambda) - \diag_i(\mu^\vee)),
\\[.15in]
\displaystyle
\Dmax(\lambda,\mu)=\diag_0(\lambda) - \max_{i=-k,\dots,n-k} 
(\diag_i(\mu^\vee) - \diag_i(S^k(\lambda))).
\end{array}
$$
The equivalence of these formulas to the definition of $\Dmin$ and $\Dmax$ 
in terms of the function $\phi_i$ is a consequence of the following identities,
which we leave as an exercise for the reader:
$$
\begin{array}{l}
\diag_{i-k}(\lambda)-\diag_{i-k}(\mu^\vee) = \phi_{i}(\mu^\vee)-
\phi_{i}(\lambda),
\\[.15in]
\phi_i(\mu^\vee)=-\phi_{-i}(\mu),
\\[.15in]
\diag_0(\lambda)=k-\phi_k(\lambda),\\[.15in]
\phi_i(S^k(\lambda)) = \phi_{i+k}(\lambda)-\phi_k(\lambda).
\end{array}
$$

\begin{proof}[Proof of Theorem~\ref{th:dmin-max}]
Let us first verify that $\Dmin\leq \Dmax$.  
We need to check that, for any integers $i$ and $j$, we have
$-\phi_i(\lambda)-\phi_{-i}(\mu)\leq 
-\phi_j(\lambda)-\phi_{-j+k-n}(\mu)$,
or, equivalently,
$$
\phi_j(\lambda)-\phi_i(\lambda)\leq \phi_{-i}(\mu)-\phi_{-j+k-n}(\mu)
=\phi_{-i}(\mu)-\phi_{-j+k}(\mu)+k.
$$
We may assume that $j\in[i,i+n[$ because the function $\phi_j$ satisfies the 
condition $\phi_{j+n}=\phi_j + k$.
Then we have $\phi_j(\lambda)-\phi_i(\lambda)\leq \min(j-i,k)$.
Indeed, $\phi_j(\lambda)-\phi_i(\lambda)\leq \phi_{i+n}-\phi_i=k$
and $\phi_j(\lambda)-\phi_i(\lambda)\leq j-i$ because 
$\phi_{s+1}-\phi_s\in\{0,1\}$ for any $s$.
On the other hand, we have
$\phi_{-i}(\mu)-\phi_{-j+k}(\mu)+k\geq \min(j-i,k)$ or, equivalently,
$ \phi_{k-j}(\mu)- \phi_{-i}(\mu) 
\leq \max(i+k-j,0)$.
Indeed, if $k-j\leq -i$ then the left-hand side is non-positive and the 
right-hand side is zero; otherwise 
$ \phi_{k-j}(\mu) - \phi_{-i}(\mu)  \leq (k-j)-(-i) = i+k-j$.
This proves the required inequality.

Let us now show that the values of $d$, for which $q^d$ occurs 
with nonzero coefficient in
$\sigma_\lambda*\sigma_\mu$, form the interval $[\Dmin,\Dmax]$.
According to Corollary~\ref{cor:main-variant},  the power $q^d$ appears
in the quantum product  $\sigma_\lambda*\sigma_\mu$ whenever
$\mu^\vee/d/\lambda$ is a valid toric shape.
This is true if and only if the following two conditions are satisfied:
(a) $\mu^\vee[d]\geq \lambda[0]$, i.e., $\mu^\vee[d]_i\geq \lambda[0]_i$ 
for all $i$; and
(b) $\lambda^\da[0^\da]\geq \mu^\vee[d]$,
where $\lambda^\da[0^\da]=S^k(\lambda)[\diag_0(\lambda)]$,
cf.~Lemma~\ref{lem:loops-torus}.
The first condition~(a) can be written as $\phi_i(\lambda)-\phi_i(\mu^\vee)+d=
\phi_i(\lambda)+\phi_{-i}(\mu)+d\geq 0$ for all $i$.
It is equivalent to the inequality $d\geq \Dmin$.
The second condition~(b) can be written as 
$\phi_i(\mu^\vee)-\phi_i(\lambda^\da)+0^\da-d=
-\phi_{-i}(\mu)-(\phi_{i+k}(\lambda)-\phi_k(\lambda))+(k-\phi_k(\lambda))-d
= -\phi_{i+k-n}(\lambda)-\phi_{-i}(\mu)-d\geq 0$ for all $i$.
It is equivalent to the inequality $d\leq \Dmax$.
\end{proof}

Recall that in Section~\ref{ssec:essential}, for a triple
of partitions $\lambda,\mu,\nu\in P_{kn}$, we defined
the essential interval $[\dmin,\dmax]$.

\begin{corollary}
For a pair of partitions $\lambda,\mu\in P_{kn}$, we have
$$
[\Dmin(\lambda,\mu),\Dmax(\lambda,\mu)]=\bigcup_{\nu\in P_{kn}}
[\dmin(\lambda,\mu,\nu),\dmax(\lambda,\mu,\nu)].
$$
\label{cor:Dmin-dmin}
\end{corollary}

\begin{proof} 
It is clear from the definitions that, for any $\lambda,\mu,\nu$,
$$
[\dmin(\lambda,\mu,\nu),\dmax(\lambda,\mu,\nu)]
\subseteq [\Dmin(\lambda,\mu),\Dmax(\lambda,\mu)].
$$
Thus the right-hand side of the formula in Corollary~\ref{cor:Dmin-dmin}
is contained in the left-hand side.  On the other hand, by 
Proposition~\ref{prop:essential}, the right-hand side contains the set of
all $d$'s such that $q^d$ appears in $\sigma_\lambda*\sigma_\mu$, which
is equal to the left-hand side of the expression, by Theorem~\ref{th:dmin-max}. 
\end{proof}

The numbers $\Dmin$ and $\Dmax$ have very simple geometric meanings
in terms of loops on the torus $\T_{kn}$.
The number $\Dmin$ is the minimal possible $d$ such that 
$\mu^\vee[d]\geq \lambda[0]$. 
In other words, 
the loop $\mu^\vee[\Dmin]$ touches (but does not cross) the 
{\it South-East\/} side of the loop $\lambda[0]$.
Also, the number $\Dmax$ is the maximal possible $d$ such that 
$\lambda^\da[0^\da]\geq \mu^\vee[d]$,
cf.~Lemma~\ref{lem:loops-torus}.
This means that the loop $\mu^\vee[\Dmax]$ 
touches the {\it North-West\/} side of the loop $\lambda[0]$.

\begin{figure}[ht]
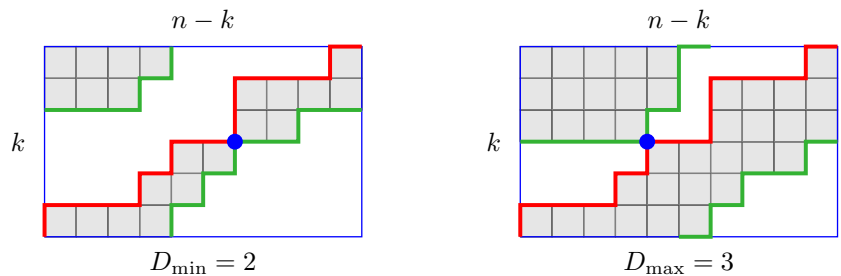

\pspicture(0,-10)(120,82)
\rput(-10,36){$k$}
\rput(60,82){$n-k$}
\rput(60,-10){$\Dmin=2$}
\ggbox(0,0)
\ggbox(1,0)
\ggbox(2,0)
\ggbox(3,0)
\ggbox(3,1)
\ggbox(4,1)
\ggbox(4,2)
\ggbox(5,2)
\ggbox(6,3)
\ggbox(7,3)
\ggbox(6,4)
\ggbox(7,4)
\ggbox(8,4)
\ggbox(9,4)
\ggbox(9,5)
\ggbox(0,4)
\ggbox(1,4)
\ggbox(2,4)
\ggbox(0,5)
\ggbox(1,5)
\ggbox(2,5)
\ggbox(3,5)
\psline[linecolor=blue, linewidth=0.5pt]{-}(0,0)(120,0)(120,72)(0,72)(0,0)
\psline[linecolor=red, linewidth=1.5pt]{-}(0,0)(0,12)(36,12)(36,24)(48,24)(48,36)(72,36)(72,60)(108,60)(108,72)(120,72)
\psline[linecolor=mygreen, linewidth=1.5pt]{-}(48,0)(48,12)(60,12)(60,24)(72,24)(72,36)(96,36)(96,48)(120,48)
\psline[linecolor=mygreen, linewidth=1.5pt]{-}(0,48)(36,48)(36,60)(48,60)(48,72)
\pscircle*[linecolor=blue](72,36){3}
\endpspicture
\qquad\qquad\qquad
\pspicture(0,-10)(120,82)
\rput(-10,36){$k$}
\rput(60,82){$n-k$}
\rput(60,-10){$\Dmax=3$}
\ggbox(0,0)
\ggbox(1,0)
\ggbox(2,0)
\ggbox(3,0)
\ggbox(4,0)
\ggbox(5,0)
\ggbox(3,1)
\ggbox(4,1)
\ggbox(5,1)
\ggbox(6,1)
\ggbox(4,2)
\ggbox(5,2)
\ggbox(6,2)
\ggbox(7,2)
\ggbox(8,2)
\ggbox(6,3)
\ggbox(7,3)
\ggbox(8,3)
\ggbox(9,3)
\ggbox(6,4)
\ggbox(7,4)
\ggbox(8,4)
\ggbox(9,4)
\ggbox(9,5)
\ggbox(0,3)
\ggbox(1,3)
\ggbox(2,3)
\ggbox(3,3)
\ggbox(0,4)
\ggbox(1,4)
\ggbox(2,4)
\ggbox(3,4)
\ggbox(4,4)
\ggbox(0,5)
\ggbox(1,5)
\ggbox(2,5)
\ggbox(3,5)
\ggbox(4,5)
\psline[linecolor=blue, linewidth=0.5pt]{-}(0,0)(120,0)(120,72)(0,72)(0,0)
\psline[linecolor=red, linewidth=1.5pt]{-}(0,0)(0,12)(36,12)(36,24)(48,24)(48,36)(72,36)(72,60)(108,60)(108,72)(120,72)

\psline[linecolor=mygreen, linewidth=1.5pt]{-}(60,00)(72,00)(72,12)(84,12)(84,24)(108,24)(108,36)(120,36)

\psline[linecolor=mygreen, linewidth=1.5pt]{-}(0,36)(48,36)(48,48)(60,48)(60,72)(72,72)
\pscircle*[linecolor=blue](48,36){3}
\endpspicture

\caption{The lowest power $\Dmin$ and the highest power $\Dmax$}
\label{fig:d-min-max}
\end{figure}

The constructions of $\Dmin$ and $\Dmax$ become transparent once 
everything is drawn on a picture.  
Figure~\ref{fig:d-min-max} gives an example for $k=6$ and $n=16$.
Here $\lambda=(9,6,6,4,3,0)$ (shown in red color)
and $\mu^\vee=(6,4,3,2,2,1)$ (shown in green color).
We have $\Dmin=2$ and $\Dmax=3$.  The left picture
shows that the green loop $\mu^\vee[\Dmin]$ touches the South-East side
of the red loop $\lambda[0]$; and the right picture shows that
the green loop $\mu^\vee[\Dmax]$ touches the North-West side of the red 
loop~$\lambda[0]$.

Let us show that the strange duality 
flips the interval $[\Dmin,\Dmax]$.
Indeed, it follows from Theorem~\ref{thm:strange-duality} that
$$
C_{\lambda\mu}^{\nu,d}=
C_{S^{n-k}(\lambda^\vee)\,\mu^\vee}^{\nu^\vee,\,\diag_0(\lambda)-d}.
$$
In other words, the coefficient of $q^d\,\sigma_\nu$ in 
the quantum product $\sigma_\lambda*\sigma_\mu$ is exactly the same 
as the coefficient of $q^{\diag_0(\lambda)-d}\,\sigma_{\nu^\vee}$
in the quantum product $\sigma_{S^{n-k}(\lambda^\vee)}*\sigma_{\mu^\vee}$.
This means that the set of all powers of $q$ that occur in 
$\sigma_\lambda*\sigma_\mu$ is obtained from the set of all powers
powers of $q$ that occur in $\sigma_{S^{n-k}(\lambda^\vee)}*\sigma_{\mu^\vee}$
by the transformation $d\mapsto \diag_0(\lambda)-d$.
In particular, we obtain the following statement.

\begin{corollary}
For any $\lambda,\mu\in P_{kn}$, we have
$$
\begin{array}{l}
\Dmin(\lambda,\mu)=\diag_0(\lambda)-\Dmax(S^{n-k}(\lambda^\vee), \mu^\vee),\\[.1in]
\Dmax(\lambda,\mu)=\diag_0(\lambda)-\Dmin(S^{n-k}(\lambda^\vee), \mu^\vee).
\end{array}
$$
\label{cor:min=max}
\end{corollary}

Recall that the map $\lambda\mapsto \tilde\lambda=S^{n-k}(\lambda^\vee))$
is an involution on $P_{kn}$ such that 
$\diag_0(\lambda)=\diag_0(\tilde\lambda)$, see the paragraph
after Theorem~\ref{th:homomorph}.

The lowest $\Dmin$ and the highest $\Dmax$ powers of $q$ in 
the quantum product $\sigma_\lambda*\sigma_\mu$ can be easily recovered
from the hidden symmetry and the strange duality of the 
Gromov-Witten invariants.
Moreover, the Gromov-Witten invariants $C_{\lambda\mu}^{\nu,d}$
in the case when $d=\Dmin$ or $d=\Dmax$ are equal to certain 
Littlewood-Richardson coefficients.

\begin{corollary}
Let $\lambda,\mu,\nu\in P_{kn}$ be three partitions.
Let $\Dmin=\Dmin(\lambda,\mu)$ and $\Dmax=\Dmax(\lambda,\mu)$.
By the definition, there are integers $a$ and $b$ such that 
$\Dmin+\phi_a(\lambda)+\phi_{-a}(\mu)=0$ and
$\Dmax+\phi_{-b}(\lambda) + \phi_{b+k-n}(\mu)=0$.   For such
$a$ and $b$, we have
$$
C_{\lambda\mu}^{\nu,\,\Dmin} = c_{S^{a}(\lambda)\,S^{-a}(\mu)}^\nu\,
\quad\textrm{and}\quad
C_{\lambda\mu}^{\nu,\,\Dmax} = \
c_{S^{b}(\lambda^\vee)\,S^{n-k-b}(\mu^\vee)}^{\nu^\vee}\,.
$$
\end{corollary}

\begin{proof}
If $\Dmin(\lambda,\mu)=\dmin(\lambda,\mu,\nu)$ then 
the statement about $C_{\lambda\mu}^{\nu,\,\Dmin}$
 is a special case of Proposition~\ref{prop:essential}.
If $\Dmin(\lambda,\mu)<\dmin(\lambda,\mu,\nu)$ then, by the same 
proposition, both sides are equal to $0$.
Similarly, the statement about $C_{\lambda\mu}^{\nu,\,\Dmax}$ follows
from Proposition~\ref{prop:essential}.
\end{proof}

This statement means that, for a {\it toric\/} shape 
$\kappa=\mu^\vee/d/\lambda$ with $d=\Dmin$,
there always exists a cyclic shift 
$S^a(\kappa)$ that is equal to the {\it skew\/} shape 
$S^a(\kappa)=S^{a}(\mu^\vee)/0/S^a(\lambda)$, cf.~Figure~\ref{fig:d-min-max}.
If $d=\Dmax$ then the same is true for the complement toric shape 
$\kappa^\vee$.

\section{Affine nil-Temperley-Lieb algebra}
\label{sec:nTL}

In this section we discuss the affine nil-Temperley-Lieb algebra
and its action on the quantum cohomology $\QH^*(Gr_{kn})$.
This section justifies the word ``affine'' 
that appeared in the title of this article.
The affine nil-Temperley-Lieb algebra presents a model for the
quantum cohomology of the Grassmannian.

\medskip

For $n\geq 2$, let us define the {\it affine nil-Temperley-Lieb algebra\/} 
$\A_n$ as the associative algebra with 1 over $\Z$ with generators
$a_i$, $i\in\Z/n\Z$, and the following defining relations:
\begin{equation}
\begin{array}{l}
a_i\, a_i =a_{i}\, a_{i+1}\, a_i = a_{i+1}\, a_i\, a_{i+1} = 0, \\ 
a_i\, a_j = a_j\, a_i,\quad\textrm{if }i-j\not\equiv \pm1.
\end{array}
\label{eq:A_n}
\end{equation}
The subalgebra of $\A_n$ generated by $a_1,\dots,a_{n-1}$ 
is called the {\it nil-Temperley-Lieb algebra}.
Its dimension is equal to the $n$-th Catalan number.
According to Fomin and Green~\cite{FG},  this algebra 
can also be defined as the algebra of operators acting on the space 
of formal combinations of Young diagrams by adding boxes to diagonals.
In the next paragraph we extend this action
to the affine nil-Temperley-Lieb algebra.

Recall $\omega(\lambda)=(\omega_1,\dots,\omega_n)$ denotes
the 01-word of a partition $\lambda\in P_{kn}$, see 
Section~\ref{sec:prelims}.  
Let us define $\lambda(\omega)\in P_{kn}$ 
as the partition with $\omega(\lambda)=\omega$.
Let $\epsilon_i$ be the $i$-th coordinate $n$-vector; and let
$\epsilon_{ij}=\epsilon_i-\epsilon_j$.
For $i,j\in\{1,\dots,n\}$, we define the $\Z[q]$-linear operator
$E_{ij}$ on the space $\QH^*(Gr_{kn})$ given in the basis of 
Schubert cells by
$$
E_{ij}:\sigma_{\lambda(\omega)}\longmapsto
\left\{
\begin{array}{cl}
\sigma_{\lambda(\omega-\epsilon_{ij})} &
\textrm{if } \omega-\epsilon_{ij} \textrm{ is a 01-word},\\[.05in]
0 & \textrm{otherwise},
\end{array}
\right.
$$
where $\omega-\epsilon_{ij}$ means the coordinatewise difference of
two $n$-vectors.
We define the action of the generators 
$a_1,\dots,a_n$ of the affine nil-Temperley-Lieb algebra $\A_n$  
on the quantum cohomology $\QH^*(Gr_{kn})$ using operators
$E_{ij}$ as follows:
$$
\begin{array}{l}
a_i = E_{i\,i+1}\quad\textrm{for } i = 1,\dots,n-1; \\[.05in]
a_n = q\cdot E_{n1}.
\end{array}
$$
It is an easy exercise to check that these operators 
satisfy relations~(\ref{eq:A_n}).

This action can also be interpreted in terms of Young diagrams
that fit inside the $k\times (n-k)$-rectangle.
For $i=1,\dots,n-1$, we have $a_i(\sigma_\lambda)=\sigma_\mu$
if the shape $\mu$ is obtained by adding a 
box to the $(i-k)$-th diagonal of the shape $\lambda$; or
$a_i(\sigma_\lambda)=0$ if it is not possible to add such a box.
Also, $a_n(\sigma_\lambda)=q\cdot \sigma_\mu$ if the shape $\mu$ is 
obtained from the shape $\lambda$ by removing a rim hook of size $n-1$;
or $a_n(\sigma_\lambda)=0$ if it is not possible to remove such a rim hook.
Notice that the partition $\mu\in P_{kn}$ is obtained from $\lambda\in P_{kn}$ 
by removing a rim hook of size $n-1$ if and only if the order ideal
$D_{\mu[r+1]}$ in the cylinder $\Cyl_{kn}$ is
obtained from $D_{\lambda[r]}$ by adding a box to the $(n-k)$-th diagonal.
Thus the generators $a_i$, $i=1,\dots,n$, of the affine nil-Temperley-Lieb 
algebra naturally act on order ideals in $\Cyl_{kn}$ by adding 
boxes to $(i-k)$-th diagonals.

\medskip
Let us say a few words on a relation between the affine 
nil-Temperley-Lieb algebra $\A_n$ and the 
{\it affine Lie algebra\/} $\slhat_n$ (without central extension).
The vector space $\H^*(Gr_{kn})\otimes\C$ 
can be regarded as the $k$-th {\it fundamental 
representation\/} $\Phi_k$ of the Lie algebra
$\sln_n$.  A Schubert class $\sigma_{\lambda(\omega)}$
corresponds to the weight vector of weight $\omega$.
These are exactly the weights obtained by conjugations of the $k$-th
fundamental weight.
The generator $e_i$ of $\sln_n$ acts on $\H^*(Gr_{kn})$ as the operator 
$a_i$ above by adding a box to the $(i-k)$-th diagonal of the shape $\lambda$.
(The generators $e_i$ of $\sln_n$ should not be confused with 
elementary symmetric functions.)
The generator $f_i$ acts as the conjugate to $e_i$ operator by removing a box 
from the  $(i-k)$-th diagonal.
Recall that every representation $\Gamma$ of $\sln_n$ 
gives rise to the {\it evaluation module\/} $\Gamma(q)$, which is
a representation of the affine Lie algebra $\slhat_n$, see~\cite{Kac}.
Then the space $\QH^*(Gr_{kn})\otimes \C[q,q^{-1}]$ can be regarded
as the evaluation module $\Phi_k(q)$ of the $k$-fundamental representation:
$$
\QH^*(Gr_{kn})\otimes \C[q,q^{-1}] \simeq \Phi_k(q).
$$
This equality is just a formal identification of two linear spaces
over $\C[q,q^{-1}]$ given by mapping a Schubert class to the 
corresponding weight vector in $\Phi_k(q)$.
This $\C[q,q^{-1}]$-linear action of $\slhat_n$ on 
$\QH^*(Gr_{kn})\otimes \C[q,q^{-1}]$ is explicitly given by 
$$
\begin{array}{l}
e_i = E_{i\,i+1}\quad\textrm{for } i = 1,\dots,n-1, \textrm{ \ and \ }
e_n = q\cdot E_{n1}, \\[.05in]
f_i = E_{i+1\,i}\quad\textrm{for } i = 1,\dots,n-1, \textrm{ \ and \ }
f_n = q^{-1}\cdot E_{1n}, \\[.05in]
h_i: \sigma_\lambda\mapsto (\omega_i(\lambda)- \omega_{i+1}(\lambda))
\sigma_\lambda \quad\textrm{for } i = 1,\dots,n, 
\end{array}
$$
where we assume that $\omega_{n+1}(\lambda)=\omega_1(\lambda)$.

Let $\n$ be the subalgebra of the affine algebra $\slhat_n$ generated 
by $e_1,\dots,e_n$.  The affine 
nil-Temperley-Lieb algebra (with complex coefficients)
is exactly the following quotient of the 
universal enveloping algebra $U(\n)$ of $\n$:
$$
\A_n\otimes\C \simeq U(\n)/\<(e_i)^2\mid i=1,\dots,n\>.
$$
Indeed, Serre's relations modulo the ideal $\<(e_i)^2\>$ 
degenerate to the defining relations~(\ref{eq:A_n}) of $\A_n$.
Notice that the squares of the generators $(e_i)^2$ and $(f_i)^2$ vanish
in all fundamental representations $\Phi_k$ and in their evaluation 
modules $\Phi_k(q)$.  
The action of the affine nil-Temperley-Lieb algebra $\A_n$ on 
$\QH^*(Gr_{kn})$ described above in this section is exactly 
the action deduced from the evaluation module $\Phi_k(q)$.

\medskip
Let us show how the affine nil-Temperley-Lieb algebra $\A_n$
is related to cylindric shapes.  Let $\kappa$ be a cylindric shape
of type $(k,n)$ for {\it some\/} $k$.
Let us pick any cylindric tableau $T$ of shape $\kappa$
and {\it standard\/} weight $\beta=(1,\dots,1)$.
For $i=1,\dots,|\kappa|$, let $d_i$ be $k$ plus the index of the diagonal 
that contains the entry $i$ in the tableau~$T$.
Let us define $\a_\kappa=a_{d_1}\cdots a_{d_{|\kappa|}}$.  
The monomials for different
tableaux of the same shape can be related by the commuting relations
$a_i\cdot a_j=a_j\cdot a_i$.
Thus the monomial $\a_\kappa$ does not depend on the choice of tableau.
For two cylindric shapes $\kappa$ and $\tilde\kappa$ of types $(k,n)$ and 
$(\tilde k,n)$, let us write $\kappa\sim \tilde\kappa$ 
whenever $\a_\kappa=\a_{\tilde\kappa}$.
Clearly, $\a_\kappa$ does not change if we shift the shape $\kappa$
in the South-East direction.  Thus $\kappa\sim\tilde\kappa$ for any 
$\tilde\kappa$ obtained from $\kappa$ by such a shift.  
Moreover, if the diagram $D_\kappa$ of $\kappa$ has several connected 
components then we can shift each connected component independently. 
These shifts of connected components generate the equivalence
relation `$\sim$'.
Any nonvanishing monomial in $\A_n$ is equal to $\a_\kappa$
for some $\kappa$.  
Thus the map $\kappa\mapsto \a_\kappa$ gives 
a one-to-one correspondence between cylindric shapes 
(modulo the `$\sim$'-equivalence) and nonvanishing monomials 
in the algebra~$\A_n$.  

For any $\mu\in P_{kn}$ and a cylindric shape $\kappa$ there is at most one 
cylindric loop $\lambda[d]$ of type $(k,n)$
such that $\lambda/d/\mu\sim\kappa$.
The action of a monomial $\a_\kappa$ 
on $\QH^*(Gr_{kn})$ is given by 
\begin{equation}
\a_\kappa:\sigma_\mu \longmapsto 
\left\{
\begin{array}{cl}
q^d\, \sigma_\lambda & \textrm{if } 
\lambda/d/\mu\sim \kappa.\\[.15in]
0 & \textrm{if there are no such $\lambda$ and $d$.}
\end{array}
\right.
\label{eq:a-kappa}
\end{equation}

\medskip
So far in this section we treated the quantum cohomology $\QH^*(Gr_{kn})$
as a linear space.  Let us show that the action
of the affine nil-Temperley-Lieb algebra $\A_n$ is helpful for
describing the multiplicative structure of $\QH^*(Gr_{kn})$.

Let us define the elements $\e_1,\dots,\e_{n-1}$ 
and $\h_1,\dots,\h_{n-1}$ in the algebra $\A_n$ 
as follows.
For a proper subset $I$ in $\Z/n\Z$, let 
$\prod_{i\in I}^\circlearrowright  a_i\in\A_n$ 
be the product
of $a_i$, $i\in I$, taken in an order such that if $i,i+1\in I$ then 
$a_{i+1}$ goes before $a_{i}$.  
This product is well-defined because all such orderings 
of $a_i$, $i\in I$, are obtained from each other by switching
commuting generators.
Also, let $\prod_{i\in I}^\circlearrowleft  a_i\in\A_n$ 
be the element obtained by reversing 
the ``cyclic order'' of $a_i$'s
in $\prod_{i\in I}^\circlearrowright  a_i$.
Let us define
$$
\e_r = \sum_{|I|=r}\ \prod_{i\in I}^\circlearrowright  a_i
\qquad\textrm{and}\qquad
\h_r = \sum_{|I|=r}\ \prod_{i\in I}^\circlearrowleft a_i\,,
$$
where the sum is over all $r$-element subsets $I$ in $\Z/n\Z$.
For example,
$$
\begin{array}{l}
\e_1 = \h_1 = a_1+\cdots+a_n,\\[.1in]
\e_2 = a_2\, a_1 + a_3\, a_2 + \cdots + a_n\, a_{n-1} + a_1\, a_n + 
\sum^{\,\mathrm{c}} a_i\, a_j,
\\[.1in]
\h_2 = a_1\, a_2 + a_2\, a_3 + \cdots + a_{n-1}\, a_n + a_n\, a_1 + 
\sum^{\,\mathrm{c}} a_i\, a_j,
\end{array}
$$
where $\sum^{\,\mathrm{c}} a_i\, a_j$
is the sum of products of (unordered) pairs of commuting $a_i$ and $a_j$,
i.e., $i$ and $j$ are not adjacent elements in $\Z/n\Z$.
In the spirit of~\cite{FG}, we can say that the $\e_r$ are
elementary symmetric polynomials
and the $\h_r$ are the complete homogeneous symmetric polynomials
in {\it noncommutative\/} variables $a_1,\dots,a_n$.
Notice that the element $\e_r$ (respectively, $\h_r$) in $\A_n$
is the sum of monomials $\a_\kappa$ for all non-`$\sim$'-equivalent 
cylindric vertical (respectively, horizontal) $r$-strips $\kappa$.

\begin{lemma}
The elements $\e_i$ and $\h_j$ in the algebra $\A_{n}$ 
commute pairwise.  For $i+j>n$, we have $\e_i\cdot \h_j=0$.
These elements are related by the equation 
\begin{equation}
\left(1+\sum_{i=1}^{n-1} \e_i\, t^i\right)\cdot
\left(1+\sum_{j=1}^{n-1} \h_j\, (-t)^j\right)=
1+\left(\sum_{k=1}^{n-1} (-1)^{n-k}\,\e_k\cdot \h_{n-k}\right) t^n.
\label{eq:e-h-AN}
\end{equation}
\label{lem:e-h-commute-in-AN}
\end{lemma}

\begin{proof}
Let us first show that $\e_i\cdot \h_j=0$, 
for $i+j>n$.  Indeed, by the pigeonhole principle, every monomial
in the expansion of $\e_i\cdot \h_j$ contains two repeating generators
$a_s$. If there is such a monomial that does not vanish in $\A_n$
then it is of the form $\a_\kappa$ and the shape $\kappa$
contains at least two elements in the $(s-k)$-th diagonal.
Thus $\kappa$ should contain a $2\times 2$ rectangle.
But it is impossible to cover a $2\times 2$ rectangle by a 
horizontal and a vertical strip.

Two elements $\h_i$ and $\h_j$ commute because the coefficient
of a monomial $\a_\kappa$ in $\h_i \cdot \h_j$ is equal to the number
of cylindric tableaux of shape $\kappa$ and weight $(i,j)$, which is 
the same as the number of tableaux of weight $(j,i)$, by 
Corollary~\ref{cor:K-Sym}.

Let us check that the coefficient of $t^l$, $0<l<n$, in the left-hand
side of~(\ref{eq:e-h-AN}) is zero.  Indeed, every nonvanishing 
monomial in the expansion of $\e_i\cdot \h_j$, $i+j<n$, correspond
to a cylindric shape $\kappa$ with $|\kappa|<n$.  Thus there is 
a cyclic shift of $\kappa$ that becomes a {\it skew\/} shape 
$\lambda/\mu$.  The formula follows from the classical result 
about adding horizontal and vertical strips 
in the planar (non-affine) case.

Finally, the relation~(\ref{eq:e-h-AN}) allows one to express the 
elements $\e_1,\dots,\e_{n-1}$ in terms of $\h_1,\dots,\h_{n-1}$,
which shows that the elements $\e_i$ commute with each other and with 
the elements $\h_j$.
\end{proof}

Recall that the quantum cohomology ring $\QH^*(Gr_{kn})$ is
the quotient (\ref{eq:QH}) of the polynomial ring over $\Z[q]$
in the variables $e_1,\dots,e_{k},h_1,\dots,h_{n-k}$. 
These generators
are the special Schubert classes 
$e_i=\sigma_{1^i}$ and $h_j = \sigma_j$\,.
We can reformulate Bertram's quantum Pieri formula, 
see Proposition~\ref{prop:Pieri}, as follows.

\begin{corollary} {\rm (Quantum Pieri formula: $\A_n$-version)}
For any $\lambda\in P_{kn}$, the products the Schubert class  
$\sigma_\lambda$ in the quantum cohomology ring $\QH^*(Gr_{kn})$ 
with the generators $e_i$ and $h_j$ are given by 
$$
e_i*\sigma_\lambda = \e_i(\sigma_\lambda)\quad
\textrm{and}\quad
h_j*\sigma_\lambda = \h_j(\sigma_\lambda),
$$
where and $i=1,\dots,k$ and $j=1,\dots,n-k$.
\label{cor:pieri-AN}
\end{corollary}

Indeed, by~(\ref{eq:a-kappa}) the operators $\e_i$ and $\h_j$
act on $\QH^*(Gr_{kn})$ by adding cylindric vertical $i$-strips and 
horizontal $j$-strips, respectively.

The quantum Giambelli formula~(\ref{eq:Giambelli})
implies the following statement.

\begin{corollary}
For any $\lambda\in P_{kn}$ the element
$$
\mathbf{s}_\lambda=\det(\h_{\lambda_i+j-i})_{1\leq i,j\leq k} =
\det(\e_{\lambda'_i+j-i})_{1\leq i,j\leq n-k} 
\in \A_n
$$
acts on the quantum cohomology $\QH^*(Gr_{kn})$ as the operator of 
quantum multiplication by the Schubert class~$\sigma_\lambda$.
The element $\mathbf{s}_\lambda\in \A_n$ is given by
the following positive linear combination of monomials:
$$
\mathbf{s}_\lambda= \sum_{\kappa=\nu/d/\mu}
C_{\lambda\mu}^{\nu,d}\, \a_\kappa\,
$$
where the sum over all non-`$\sim$'-equivalent cylindric shapes $\kappa$.
\label{cor:det}
\end{corollary}

The second claim follows from~(\ref{eq:a-kappa}).
Thus, even though the expansion of the determinant contains negative
signs, all negative terms cancel, and 
$\mathbf{s}_\lambda$ always reduces to a {\it positive\/} expression.

\medskip
The algebra $\A_n$ acts on $\QH^*(Gr_{kn})$ for all values
of $k$.  In order to single out one particular $k$, we need to 
describe certain $n-1$ central elements in the algebra $\A_n$.
We say that a cylindric shape $\kappa$ of type $(k,n)$
is a {\it circular ribbon\/} if the diagram of $\kappa$
contains no $2\times 2$ rectangle and $|\kappa|=n$.
Up to the `$\sim$'-equivalence, 
there are exactly $\binom nk$ circular ribbons of type $(k,n)$.
Let us define the elements $\zz_1,\dots,\zz_{n-1}$ in 
$\A_n$ as the sums
$$
\zz_k = \sum_{\kappa} \a_\kappa
$$
over all $\binom nk$ non-`$\sim$'-equivalent circular
ribbons $\kappa$ of type $(k,n)$.
These elements are also given by
$$
\zz_k = \e_k\cdot \h_{n-k}\,.
$$
Indeed, a nonvanishing monomial in $\e_k\cdot \h_{n-k}$ 
should be of the type $\a_\kappa$, where $\kappa$ contains
no $2\times 2$-rectangle, cf.~proof of Lemma~\ref{lem:e-h-commute-in-AN}.
Since $|\kappa|=k+(n-k)=n$, the cylindric shape $\kappa$ should be a 
circular ribbon.  Then each circular ribbon of type $(k,n)$ uniquely 
decomposes into a product of two monomials
corresponding to a vertical $k$-strip and a horizontal $(n-k)$-strip.

\begin{lemma}
The elements $\zz_1,\dots,\zz_{n-1}$ are central elements in the 
algebra $\A_n$.  For $k\ne l$, we have $\zz_k\cdot \zz_l = 0$. 
\end{lemma}

\begin{proof}
For any $i$, both elements $\zz_k\cdot a_i$ and $a_i\cdot \zz_k$ are
given by the sum of monomials $\a_\kappa$ over all cylindric shapes $\kappa$,
$|\kappa|=n+1$, that have exactly one $2\times 2$ rectangle centered in
the $(i-k)$-th diagonal.  Thus  $\zz_k\cdot a_i = a_i\cdot \zz_k$,
for any $i$; which 
implies that $\zz_k$ is a central element in $\A_n$.
The second claim follows from~(\ref{eq:a-kappa}).
\end{proof}

Let us define the algebra $\A_{kn}$ as
$$
\A_{kn} = \A_n\otimes\Z[q,q^{-1}]/\<\zz_1,\dots,\zz_{k-1},\zz_k-q,\zz_{k+1},
\dots, \zz_{n-1}\>. 
$$

\begin{proposition}
The localization of the
quantum cohomology $\QH^*_{\<q\>}(Gr_{kn})
=\QH^*(Gr_{kn})\otimes\Z[q,q^{-1}]$
is isomorphic to the subalgebra of $\A_{kn}$ generated by the elements 
$\e_i$ and/or $\h_j$.  This isomorphism is given by 
the $\Z[q,q^{-1}]$-linear map that sends the generators 
$e_i$ and $h_j$ of $\QH^*_{\<q\>}$ to the elements $\e_i$ and $\h_j$ 
in $\A_{kn}$, respectively.
\end{proposition}

\begin{proof}
By Corollary~\ref{cor:pieri-AN}, the algebra $\A_{kn}$ acts faithfully on 
$\QH^*_{\<q\>}(Gr_{kn})$.  The only thing that we need to check is
that the elements
$\e_i$ and $\h_j$ in $\A_{kn}$ satisfy the same relations as
the elements $e_i$ and $h_j$ in the quantum cohomology do, 
cf.~(\ref{eq:QH}).  The right-hand side of
the equation~(\ref{eq:e-h-AN}) becomes $1+(-1)^{n-k} q\, t^n$
in the algebra $\A_{kn}$. 
It remains to show that $\e_i=\h_j=0$ in $\A_{kn}$ whenever
$i>k$ and $j>n-k$.  By Lemma~\ref{lem:e-h-commute-in-AN}, we have 
$\e_i\cdot \h_{n-k}=\h_j\cdot \e_k=0$, for $i>k$ and $j>n-k$.
Since $\zz_k=\e_k \cdot \h_{n-k}=q$, both elements $\e_k$ and 
$\h_{n-k}$ are invertible in $\A_{kn}$.  Thus $\e_i=\h_j=0$
as needed.
\end{proof}

\begin{remark} 
{\rm
Fomin and Kirillov~\cite{FK} defined a certain quadratic algebra 
and a set of its pairwise commuting elements, called {\it 
Dunkl elements}.  According to {\it quantum Monk's formula\/} from~\cite{FGP},
the multiplication in the quantum cohomology ring $\QH^*(Fl_n)$ of  the
{\it complete flag manifold\/} $Fl_n$ can be written in terms of the 
Dunkl elements.  
A conjecture from~\cite{FK}, which was proved in~\cite{P1}, says
that these elements 
generate a subalgebra isomorphic to $\QH^*(Fl_n)$.  
This section shows that the affine nil-Temperley-Lieb algebra 
$\A_n$ is, in a sense, a Grassmannian analogue of 
Fomin-Kirillov's quadratic algebra.  The pairwise commuting elements
$\e_i$ and $\h_j$ are analogues of the Dunkl elements.
It would be interesting to extend these two opposite cases
to the quantum cohomology of an arbitrary partial flag manifold.
}
\end{remark}

\section{Open questions, conjectures, and final remarks}
\label{sec:final}

\subsection{Quantum Littlewood-Richardson rule.}

The problem that still remains open is to give a generalization
of the Littlewood-Richardson rule to the quantum cohomology ring
of the Grassmannian.  As we already mentioned, it is possible to use 
the quantum Giambelli formula in order to derive a rule 
for the Gromov-Witten invariants 
$C_{\lambda\mu}^{\nu,d}$ that involves an {\it alternating\/} sum, e.g., 
see~\cite{BCF} or Corollary~\ref{cor:det} in the present paper.
The problem is to present a {\it subtraction-free\/} rule for the 
Gromov-Witten invariants.
In other words, one would like to get a combinatorial or algebraic construction
for the Gromov-Witten invariants that would imply their nonnegativity.
There are several possible approaches to this problem. 
Buch, Kresch, and Tamvakis~\cite{BKT} showed that the Gromov-Witten invariants
of Grassmannians are equal to some intersection numbers for two steps flag 
manifolds; and they conjectured a rule for the latter numbers.

In the next subsection we propose an algebraic approach to this problem
via representations of symmetric groups.

\subsection{Toric Specht modules}

For any toric shape $\kappa=\lambda/d/\mu$, let us define a
representation $S^\kappa$ of the symmetric group $S_N$, where $N=|\kappa|$,
as follows.
Let us fix a labelling of the boxes of $\kappa$ by numbers $1,\dots,N$.
Recall that every toric shape has rows and columns,
see Section~\ref{sec:cylindric-tableaux}.
The rows (columns)
of $\kappa$ give a decomposition of $\{1,\dots,N\}$ into a union of disjoint
subsets.  Let $R_\kappa\subset S_N$ and $C_\kappa\subset S_N$ 
be the {\it row stabilizer\/} and the {\it column stabilizer},
correspondingly. 
Let $\C[S_N]$ denote the group algebra of the symmetric group $S_N$.
The {\it toric Specht module\/} $S^\kappa$ is defined the 
subspace of $\C[S_N]$ given by
$$
S^\kappa=\left(\sum_{u\in R_\kappa} u \right) 
\left(\sum_{v\in C_\kappa} (-1)^{\mathrm{sign}(v)} v\right) \C[S_N].
$$
It is equipped with the action of $S_N$ by left multiplications.

If $\kappa$ is a usual shape $\lambda$ then $S^\lambda$ is known to
be an irreducible representation of $S_N$.
The following conjecture proposes how the $S_N$-module $S^{\kappa}$ decomposes
into irreducible representations, for an arbitrary toric shape $\kappa$.

\begin{conjecture}  
For a toric shape $\kappa=\lambda/d/\mu$,
the coefficients of irreducible components in the toric Specht module 
$S^{\lambda/d/\mu}$ are the Gromov-Witten invariants:
$$
S^{\lambda/d/\mu} = \bigoplus_{\nu\in P_{kn}} C_{\mu\nu}^{\lambda,d}\,S^\nu.
$$
\end{conjecture}

Equivalently, the toric Specht module $S^{\lambda/d/\mu}$
is expressed in terms of the irreducible modules $S^\nu$ in exactly the same 
way how the toric Schur polynomial $s_{\lambda/d/\mu}$ is expressed in
terms of the usual Schur polynomials $s_\nu$.

We verified this conjecture for several toric shapes.  For example,
it is easy to prove the conjecture for $k\leq 2$.
If the conjecture is true in general, it would provide an algebraic 
explanation of nonnegativity of the Gromov-Witten invariants.

Remark that Reiner and Shimozono~\cite{RS} investigated Specht modules 
for some class of shapes, called {\it percent-avoiding}, that is more 
general than skew shapes.  
A toric shape, however, may not be percent-avoiding.

\subsection{Representations of $\mathrm{GL}(k)$ and crystal bases}

According to Theorem~\ref{th:main}, each toric Schur polynomial 
$s_{\lambda/d/\mu}(x_1,\dots,x_k)$ is Schur-positive.  The usual Schur 
polynomials in $k$ variables are the characters of irreducible representations
of the general linear group $\mathrm{GL}(k)$.  Thus we obtain the following statement.

\begin{corollary}
For any toric shape $\lambda/d/\mu$, there exists a representation 
$V_{\lambda/d/\mu}$ of $\mathrm{GL}(k)$
such that $s_{\lambda/d/\mu}(x_1,\dots,x_k)$ is the character of 
$V_{\lambda/d/\mu}$.
\end{corollary}

It would be extremely interesting to present a more explicit construction 
for this representation~$V_{\lambda/d/\mu}$.  

Recall that with every representation of $\mathrm{GL}(k)$ it is possible to associate 
its {\it crystal}, which is a certain directed graph with labelled edges, e.g.,
see~\cite{KN}.  This graph encodes the corresponding representation
of $U_q(\mathfrak{gl}_k)$ modulo $\<q\>$.  Its vertices correspond to the 
elements of certain preferable basis, called the {\it crystal basis},
and the edges describe the action of generators on the basis elements.
It is well-known, e.g., see~\cite{KN}, that crystals are intimately
related to the Littlewood-Richardson rule.

The vertices of the crystal for $V_{\lambda/d/\mu}$ should correspond to 
the toric tableaux of shape $\lambda/d/\mu$.  Its edges should connect 
the vertices in a certain prescribed manner.  In a recent paper~\cite{Stem}
Stembridge described simple local conditions that would ensure that a given
graph is a crystal of some representation.
Thus in order to find the crystal for $V_{\lambda/d/\mu}$ 
it would be enough to present
a graph on the set of toric tableaux that complies with Stembridge's conditions.

Actually, an explicit construction of the crystal for 
$V_{\lambda/d/\mu}$ 
would immediately produce the following subtraction-free combinatorial
rule for the Gromov-Witten invariants:  The Gromov-Witten invariant 
$C_{\mu\nu}^{\lambda,d}$ is equal to the number of toric tableaux $T$ of
shape $\lambda/d/\mu$ and weight $\nu$ such that there are no 
directed edges in the crystal with initial vertex $T$.
The last condition means that the element in the crystal basis
given by $T$ is annihilated by the operators $\tilde{e}_i$.

\medskip
Remark that all numerous (re)formulations of the Littlewood-Richardson
rule and all explicit constructions of crystals for representations
of $\mathrm{GL}(k)$ use some kind of ordering of elements in shapes.
The main difficulty with toric shapes is that they are cyclically 
ordered and there is no natural way to select a linear order on a cycle.



\subsection{Verlinde algebra and fusion product}

Several people observed that the specialization
of the quantum cohomology ring $\QH^*(Gr_{kn})$ at $q=1$ is isomorphic
to the {\it Verlinde algebra\/} (a.k.a.\ the {\it fusion ring})
of $\mathrm{U}(k)$ at level $n-k$, see
Witten~\cite{Wit} for a physical proof and Agnihotri~\cite{Agni} for
a mathematical proof.  This ring is the Grothedieck ring of representations 
of $\mathrm{U}(k)$ modulo some identifications.
A Schubert class $\sigma_\lambda$ corresponds to the irreducible
representation $V_\lambda$ with highest weight given by 
the partition $\lambda$.

All constructions of this article for the quantum product make perfect 
sense for the Verlinde algebra and its product, called the 
{\it fusion product.}
Our strange duality might have a natural explanation in terms of the
Verlinde algebra.

\subsection{Geometrical interpretation}

The relevance of skew Young diagrams to the product of Schubert
classes in the cohomology ring $\H^*(Gr_{kn})$ has a geometric explanation, 
see~\cite{Fulton}.  It is possible to see that the intersection
of two Schubert varieties $\Omega_\lambda\cap \tilde\Omega_\mu$ 
(where $\tilde\Omega_\mu$ is taken in the opposite Schubert 
decomposition)
is empty unless $\mu/\lambda^\vee$ is a valid skew shape.
A natural question to ask is:  How to extend this construction
to the quantum cohomology ring $\QH^*(Gr_{kn})$ and toric shapes?
It would be interesting to obtain a ``geometric''  proof of 
our result on toric shapes (Corollary~\ref{cor:main-variant}), 
and also to present a geometric explanation of the strange duality 
(Theorem~\ref{thm:strange-duality}).

\subsection{Generalized flag manifolds}

The main theorem of~\cite{FW}  is given in a uniform setup
of the generalized flag manifold $G/P$, where $G$ 
is a complex semisimple Lie group and $P$ is its parabolic subgroup.  
It describes the minimal monomials $q^d$ in the quantum parameters $q_i$ 
that occur in the quantum product of two Schubert classes.  
It would be interesting to describe all monomials $q^d$ that occur 
with nonzero coefficients in a quantum product.

In~\cite{P3}, we proved several results for $G/B$,
where $B$ is a Borel subgroup.  We showed
that there is a unique minimal monomial $q^d$ that occurs
in a quantum product.  This monomial has a simple interpretation in terms 
of directed paths in the quantum Bruhat graph from~\cite{BFP}.
For the flag manifold $\mathrm{SL}(n)/B$, we gave a complete
characterization of all monomials $q^d$ that occur in a quantum 
product. In order to do this, we defined path Schubert polynomials
in terms of paths in the quantum Bruhat graph and showed that their
expansion coefficients in the basis of usual Schubert polynomials
are the Gromov-Witten invariants for the flag manifold.

In forthcoming publications we will address the question of
extending the constructions of~\cite{P3} and of the present paper 
to the general case~$G/P$.


\end{document}